\documentclass{imsart}

\RequirePackage{amsthm,amsmath,amsfonts,amssymb}
\RequirePackage[authoryear]{natbib}
\RequirePackage[colorlinks,citecolor=blue,urlcolor=blue,backref=page,backref=page]{hyperref}
\RequirePackage{graphicx}

\startlocaldefs
\RequirePackage{amsthm,amsmath,amsfonts,amssymb}

\usepackage{pdflscape}
\usepackage{booktabs}
\usepackage{amsmath}
\usepackage{amsthm}
\usepackage{amsfonts, mathrsfs}
\usepackage{mathtools}
\usepackage{amssymb}
\usepackage[makeroom]{cancel}
\usepackage{graphicx,caption}
\usepackage{tikz}
\usepackage{tikz-cd}
\usetikzlibrary{shapes,arrows}
\usepackage{bigints}
\usepackage{color}
\usepackage{float}

\usepackage{enumitem}
\usepackage{hyperref}
\hypersetup{bookmarksdepth=4}
\hypersetup{pdfcreator=XeLaTeX with hyperref}
\usepackage{bbm}
\hypersetup{
    colorlinks,
    citecolor=teal,
    filecolor=blue,
    linkcolor=black,
    urlcolor=teal
}
\usepackage{xcolor}
\newcommand{\norm}[1]{\left\lVert#1\right\rVert}

\numberwithin{equation}{section}
\theoremstyle{plain}
\newtheorem{theorem}{Theorem}[section]

\newtheorem{assumption}{Assumption}

\theoremstyle{remark}

\newtheorem{example}{Example}[section]

\newtheorem*{assum*}{\assumptionnumber}
\providecommand{\assumptionnumber}{}


\endlocaldefs

\begin{document}

\begin{frontmatter}
\title{A Bernstein--von Mises Theorem for Generalized Fiducial Distributions}
\runtitle{A Bernstein--von Mises Theorem for Generalized Fiducial Distributions}

\begin{aug}
\author[A]{\fnms{J.E.}~\snm{Borgert}\ead[label=e1]{elyseb@live.unc.edu}\orcid{0009-0002-8083-8358}}
\and
\author[A]{\fnms{Jan}~\snm{Hannig}\ead[label=e2]{jan.hannig@unc.edu}\orcid{0000-0002-4164-0173}}

\address[A]{Statistics and Operations Research,
University of North Carolina at Chapel Hill\printead[presep={,\ }]{e1,e2}}

\end{aug}

\begin{abstract}
An established and growing literature on generalized fiducial inference and related fiducial ideas points to the adoption of fiducial inference as a mainstream perspective among modern statisticians. Like Bayesian posteriors, generalized fiducial distributions (GFDs) are known to satisfy Bernstein–von Mises (BvM)-type results under classical regularity conditions. Existing fiducial BvM results, however, rely on relatively restrictive smoothness assumptions and are limited in scope. In this paper, we establish a Bernstein--von Mises theorem for generalized fiducial inference under the general framework of local asymptotic normality, which accommodates non-i.i.d. data settings and reduces to the familiar differentiability in quadratic mean condition in the i.i.d. case. We apply our result to extend existing fiducial theory for free-knot spline models first developed in \citet{sonderegger.2014.springer}, and further illustrate its generality in models where classical regularity conditions fail or i.i.d. assumptions are not met.
\end{abstract}



\end{frontmatter}

\section{Introduction}\label{s:intro}
Prior-driven inference is a source of long-standing debate and, in many cases, of practical challenge in the statistical community. One solution, first proposed by R.A. Fisher, is to altogether avoid the use of priors in constructing a posterior-like distribution, an approach which he termed \textit{fiducial inference} \citep{fisher.1935.eugenics}. Fisher's idea was to derive probability statements about the unknown parameters using only probabilities from the data model by switching the roles of the parameters and data; this is a similar idea to switching the roles of parameters and data to obtain a likelihood function from the data density function. Despite its initial lack of endorsement due to non-exactness issues in multivariate cases, the approach has regained favor through the development of \textit{generalized fiducial inference} (GFI) \citep{hannig.2016.tandf, murph.2023.arXiv}, which extended Fisher's fiducial idea using mechanics motivated by both generalized pivotal quantities \citep{weerahandi.1995.springer,hannig.2006.jasa} and Fraser's structural inference \citep{fraser1966structural}. 

An established and growing literature on GFI and other fiducial ideas points to the adoption of fiducial inference as a mainstream perspective among modern statisticians. \citet{liang2025extended} proposed an approach to fiducial inference that leverages deep neural networks called \textit{extended fiducial inference}, and \citet{bissiri2025new} proposed a sequential framework for generating what they call a \textit{Doob fiducial distribution}. \citet{williams2023model} showed that conformal prediction, a framework for uncertainty quantification of arbitrary point  predictions, can be derived from  a particular generalized fiducial distribution, placing GFI in the highly topical machine learning literature. That context is an interesting special case where the GFI procedure is agnostic of the underlying data model and the resulting generalized fiducial distribution (GFD) is an \textit{imprecise} probability distribution. The important implication of the model-free GFD's model agnosticism and imprecision is that it yields prediction sets that have finite-sample type I error control. A sample of other recent topics in the fiducial literature includes nonparametric methods \citep{cui2019nonparametric}, model selection \citep{koner2023eas}, decision theory \citep{taraldsen2024fiducial,williams2024decision}, precision medicine models \citep{kim2025extended}, and linear mixed effects models \citep{yang2025fiducial}. 

Generalized fiducial distributions, like Bayesian posteriors, also enjoy the \textit{Bernstein--von Mises} property. The Bernstein--von Mises (BvM) theorem asserts that under mild smoothness conditions, Bayesian posterior (or fiducial) distributions are asymptotically Gaussian, centered at locally sufficient statistics,
and with variance inversely proportional to the Fisher information (i.e. the Cram\'er-Rao lower bound). An implication of this property is that posterior distributions are asymptotically correct in a repeated sampling sense and efficient, providing frequentist justification for Bayesian (and fiducial) inference. In \citet{hannig.2009.statistica}, sufficient conditions for the univariate GFD to converge to the Bayesian posterior were given, which was in turn shown to be close to the distribution of the maximum likelihood estimator using the classical regularity conditions found in \citet{cramer.1999.book,lehmann.2006.springer}. The approach in \citet{hannig.2009.statistica} was generalized to the multivariate setting in \citet{sonderegger.2014.springer}, where it was studied in the context of free-knot spline models. 

The objective of the present work is to extend the fiducial BvM result under the general framework of local asymptotic normality (LAN). In the independent and identically distributed (i.i.d.) setting, our theorem reduces to the familiar differentiability in quadratic mean (DQM) condition established for Bayesian posteriors in \citet{lecam.1986.springer}, and later, \citet{vaart.1998.cambridge}. With the growing attention on fiducial methods, the need to establish a more general fiducial BvM result is particularly timely and, we argue, an important contribution to the literature on posterior inference. 

The result is fairly simple to show using the machinery typical of the Bayesian case; this is because the generalized fiducial distribution elicits a data-dependent prior measure whose corresponding posterior distribution locally uniformly approximates the GFD. We emphasize that this is a fully data-dependent elicitation, and that no pre-specification of a prior distribution as in the Bayesian construction is required. The remainder of the paper is organized as follows. We provide some brief background on the Bernstein--von Mises property in Section~\ref{ss:bvm} and review the mechanics of GFI in Section~\ref{ss:gfi}. Our main result is presented in Section~\ref{ss:mainresult} and is used in Section~\ref{ss:splineex} to extend the fiducial theory for free-knot spline models first developed in \citet{sonderegger.2014.springer}. The paper concludes with two illustrative numerical examples.

\section{Background}\label{s:review}

\subsection{Bernstein--von Mises results}\label{ss:bvm}
The BvM phenomenon is named for Sergei Bernstein \citep{bernstein.1917.thprb} and Richard von Mises \citep{vonmises.1931.book} and even dates back to Laplace \citep{laplace.1810.memoire}, though Doob provided the first formal proof \citep{doob.1949.application}. The classical BvM result was proved in the context of a correctly specified, i.i.d., and regular parametric model in which the dimension of the parameter is fixed and finite. Versions of the result based on the same set of regularity conditions used for establishing asymptotic normality of maximum likelihood estimators were developed in \citet{lecam.1953.uc}, \citet{walker.1969.jrssb} and \citet{dawid.1970.cambridge}, while a similar limiting theorem for posterior means was proved in \citet{bickel.1969.springer}. It was later shown in the seminal paper \citet{lecam.1970.annals} that the classical conditions requiring multiple derivatives of the log-density could be replaced by a more general condition, \textit{differentiability in quadratic mean} (DQM), which requires a quadratic approximation of the square-root densities to exist only in an average sense, an even weaker assumption than first differentiability. Le Cam's approach combines the DQM condition with a testability condition for distinguishing the true parameter akin to those established in results on posterior consistency \citep{lecam.1960.annals, schwartz.1965.springer}.  Le Cam's version of the theorem was later improved and simplified in \citet{vaart.1998.cambridge}. There, the theorem roughly states that if \(( y_1, \ldots, y_n) \) is an i.i.d. random sample from the distribution \( P_{\boldsymbol{\theta}_0}\) having density \( p_{\boldsymbol{\theta}_0} \), the model \( \boldsymbol{\theta} \mapsto p_{\boldsymbol{\theta}} \) is DQM, and the prior is continuous and places positive mass around the true parameter \( \boldsymbol{\theta}_0 \), then, 
\[
\norm{\Pi_{\bar{\boldsymbol{\Theta}} | y_1, \ldots, y_n}- N(\hat{\boldsymbol{\theta}}_n,\frac{1}{n}I_{\boldsymbol{\theta}_0}^{-1})}_{TV} \overset{P_{\boldsymbol{\theta}_0}^{n}}{\rightarrow} 0,
\]  
where \(\hat{\boldsymbol{\theta}}_n \) is any estimator satisfying \( \sqrt{n}(\hat{\boldsymbol{\theta}}_n -\boldsymbol{\theta}_0) \xrightarrow{d} N(0, I_{\boldsymbol{\theta}_0}^{-1}) \)  and \(I_{\boldsymbol{\theta}_0}\) is the Fisher information.

The literature explores asymptotic results for posterior distributions in various settings. Several papers established BvM results in the semi-parametric setting, which gave careful treatment of the complexities involving the choice of the prior and interactions between prior and likelihood. Examples of results for specific models were explored in \citet{kim.2006.annals}, \citet{de.2009.bernstein}, \citet{leahu.2011.ejs}, \citet{knapik.2011.annals}, \citet{castillo.2012.sankhya}, and \citet{kruijer.2013.ejs}, while results suitable for general models and/or general priors can be found in \citet{shen.2002.jasa}, \citet{castillo.2012.springer}, \citet{rivoirard.2012.annals}, and \citet{bickel.2012.annals}.
These developments showed that while certain priors perform well for specific functionals of interest, they may not be as effective for others. More recent work has improved upon these limitations via post-processing of the posterior, which allows for general BvM conditions that do not depend strongly on the prior and likelihood \citep{castillo.2015.semibvm}.
A very general BvM result for finite sample semi-parametric models was presented in \citet{panov.2015.bayes}, which allowed for the full parameter dimension to grow with the sample size. Representative but not exhaustive BvM results for parameters with growing dimension are \citet{ghosal.1999.bernoulli}, \citet{boucheron.2009.ejs}, and \citet{bontemps.2011.aos}. While a BvM result for infinite-dimensional parameters is not possible in standard spaces, some developments of possible notions of a non-parametric BvM property have been made \citep{castillo.2013.annals,leahu.2011.ejs}. Extensions to dependent settings are given in \citet{heyde.1979.jrss} and \citet{sweeting.1987.jrss}, while a more recent manuscript established a BvM result for a broad class of weakly dependent models \citep{connault.2021.working}. Extensions of the BvM result to misspecified models were covered in \citet{kleijn.2012.ejs} and the recent review article \citet{bochkina.2019.review}. Recent improvements on the upper bound of the dimension have substantially extended the BvM result to regimes $n >> d^2$ and demonstrated finite-sample statements \citep{katsevich.2024.arxiv}. 

\subsection{Mechanics of generalized fiducial inference}\label{ss:gfi}
GFI works by inverting a deterministic data generating algorithm (DGA) to define a data-dependent measure on the parameter space. This approach associates the data $\textbf{Y} =(Y_1,\dots,Y_n) $ with the unknown \textit{fixed} parameter $\boldsymbol{\theta}$ and an auxiliary variable $\textbf{U} = (U_1,\dots,U_n)$ through the DGA, $\textbf{Y}=G(\textbf{U}, \boldsymbol{\theta})$. The auxiliary variable \( \textbf{U} \) is assumed to have a known distribution $P_{\mathbf U}$, and can be usefully thought of as an extension of a pivotal quantity in the frequentist construction of confidence intervals and regions, where \( \textbf{Y} \) would be replaced by a sufficient statistic for \( \boldsymbol{\theta} \). In fact, this is how Fisher first approached fiducial inference. By inverting the DGA for observed, fixed data $\textbf{y} = (y_1,\dots,y_n)$, we obtain the \textit{generalized fiducial distribution} (GFD). Since the DGA is a function of a random variable with a known distribution, 
it determines a statistical model $P_{\boldsymbol{\theta}}$ for 
$\mathbf Y$ with density $p_{\boldsymbol{\theta}}(\mathbf y)$ with respect to some dominating measure.  Solving the DGA for $\boldsymbol{\theta}$ results in a distributional estimator (instead of a point or interval estimator) of the unknown parameter of interest. 

An intuitive explanation of the construction is as follows. A smoothly varying DGA $G(\textbf{U}, \boldsymbol{\theta})$ can be locally approximated as a linear function near the observed data value $\textbf{y}$. Thus, given an independent copy of $\textbf{U}$, denoted $\textbf{U}^*$, there exists a value $\boldsymbol{\theta}^*$ such that $G(\textbf{U}^*, \boldsymbol{\theta}^*)$ is closest to $\textbf{y}$. The GFD can be computed as the distribution of $\boldsymbol{\theta}^*$ using the implicit function theorem, based on the distribution of $\textbf{U}^*$ conditional on the event $\{G(\textbf{U}^*, \boldsymbol{\theta}^*) \approx \textbf{y}\}$. We emphasize that from the fiducial perspective, the parameter $\boldsymbol{\theta}$ is \textit{fixed} but unknown; the data, once observed, are also considered \textit{fixed}. The only source of randomness is the auxiliary variable \( \textbf{U} \), which is emulated in the GFD construction by an independent draw \(\textbf{U}^*\) from its known distribution. Essentially, the GFI procedure characterizes the set of $\boldsymbol{\theta}$ values that are compatible with the observed data under the data generating mechanism $G(\textbf{U}, \boldsymbol{\theta})$. This leads to a distributional estimator for $\boldsymbol{\theta}$, and while $\boldsymbol{\theta}$ is still treated as a fixed unknown, the distribution reflects uncertainty due to the randomness in \( \textbf{U} \).

To make this formal, let the data \( \textbf{y} \in \mathcal{Y} \) be observed and fixed. Define the inverse mapping
\begin{equation}\label{eq:pseudoInverse}
    Q_{\textbf{y}}(\textbf{u}) = \arg\min_{\boldsymbol{\theta}^* \in \boldsymbol{\Theta}} \norm{\textbf{y} - G(\textbf{u}, \boldsymbol{\theta}^*)},
\end{equation}
where \( \norm{\cdot} \) is typically either the \( l_2 \) or \( l_\infty \) norm. For \( \varepsilon > 0 \), consider the event
\[
\mathcal{U}_{\textbf{y},\varepsilon} = \left\{ \textbf{u} : \min_{\boldsymbol{\theta}} \norm{\textbf{y} - G(\textbf{u}, \boldsymbol{\theta})} \leq \varepsilon \right\},
\]
and let \( \textbf{U}^*_\varepsilon \) follow the distribution $P_{\mathbf U}$ truncated to $\mathcal{U}_{\textbf{y},\varepsilon}$, i.e., for a measurable $A$, the probability $P(\textbf{U}^*_\varepsilon\in A)=P(\mathbf{U}\in A\mid \textbf{U}\in\mathcal U_{\textbf{y},\varepsilon}).$ 
Denote the distribution of $Q_{\textbf{y}}(\textbf{U}^*_\varepsilon)$ by $\mu_\varepsilon$. If the \textit{weak limit} $\lim_{\varepsilon\rightarrow 0} \mu_\varepsilon$ exists, the limit is called a {\em generalized fiducial distribution.}

Note that the truncation modifies the distribution of \( \textbf{U}^*\) to only consider values for which an approximate inverse exists. This accept-reject construction of the GFD bears a strong resemblance to Approximate Bayesian Computation (ABC) methods \citep{beaumont.2002.abc}. In a typical ABC procedure, one first samples \( \boldsymbol{\theta}^* \) from the prior and generates artificial data using a data-generating algorithm \( \mathbf{y}^* = G(\mathbf{U}^*, \boldsymbol{\theta}^*) \). If the simulated data is sufficiently close to the observed data, i.e., if \( \norm{\mathbf{y} - \mathbf{y}^*} < \varepsilon \), then \( \boldsymbol{\theta}^* \) is accepted; otherwise, it is rejected. As \( \varepsilon \to 0 \), the ABC approximation converges in distribution to the true posterior distribution. The GFD procedure exhibits a similar accept-reject mechanism, but rather than sampling from a prior, one first generates \( \mathbf{U}^* \sim P_{\mathbf{U}} \), and then determines the best-fitting parameter value via
\(\boldsymbol{\theta}^* = \arg\min_{\boldsymbol{\theta}^* \in \boldsymbol{\Theta}} \norm{\mathbf{y} - G(\mathbf{U}^*, \boldsymbol{\theta}^*)}.\)
 The synthetic data \( \mathbf{y}^* = G(\mathbf{U}^*, \boldsymbol{\theta}^*) \) are then compared to the observed data, and \( \boldsymbol{\theta}^* \) is accepted if \( \norm{\mathbf{y} - \mathbf{y}^*} < \varepsilon \).

Under mild conditions (Assumptions A.1--A.4 in
 \citet{hannig.2016.tandf}), the GFD has density
\begin{equation}\label{eq:GFdensity}
    r(\boldsymbol{\theta} |\textbf{y}) = \dfrac{p_{\boldsymbol{\theta}}(\mathbf y)J(\textbf{y},\boldsymbol{\theta})}{\int_{\boldsymbol{\Theta}} p_{\boldsymbol{\theta}'}(\mathbf y)J(\textbf{y},\boldsymbol{\theta}')d\boldsymbol{\theta}'}
\end{equation}
where $J(\textbf{y},\boldsymbol{\theta}) = D\,\nabla_{\boldsymbol{\theta}}G(\textbf{u},\boldsymbol{\theta})\Big|_{\textbf{u}=G^{-1}(\textbf{y},\boldsymbol{\theta})}$. The gradient matrix $\nabla_{\boldsymbol{\theta}}G(\textbf{u},\boldsymbol{\theta})$ is computed with respect to $\boldsymbol{\theta}$, and $D$ is a determinant-like operator that depends on the norm used in Eq. \eqref{eq:pseudoInverse}. When the $l_2$ norm is used, then for an arbitrary matrix $M$, $D\,M = (\det(n^{-1}M^{\top}M))^{1/2}$.

\begin{example}\label{ex:linregression_GFsolution}(\textit{Linear Regression}). 
To demonstrate the mechanics for deriving GFDs,
consider the following simple example. The DGA for linear regression is
$$\textbf{Y} = G(\textbf{U}, \boldsymbol{\theta}) = \textbf{X} \boldsymbol\beta + \sigma\textbf{U},$$
where $\textbf{Y}$ represents the dependent variables, $\textbf{X}$ is the design matrix, $\boldsymbol{\theta} = (\boldsymbol\beta,\sigma)$ are the unknown parameters, and $\textbf{U}$ is a random variate with known density $p_{\mathbf U}$ independent of any parameters.

To compute the GFD for the linear regression model, note that the inverse map gives \( \textbf{u} = \sigma^{-1}(\textbf{y} - \textbf{X}\boldsymbol\beta) \). The gradient matrix is then
\[
\nabla_{\boldsymbol{\theta}} G(\textbf{u}, \boldsymbol{\theta}) = \left( \frac{\partial G}{\partial \boldsymbol\beta}, \frac{\partial G}{\partial \sigma} \right) = (\textbf{X}, \textbf{u}),
\]
so the Jacobian function using the \( l_2 \)-norm becomes

\[
J(\textbf{y}, \boldsymbol{\theta}) = {n^{-1}\det}\left((\textbf{X}, \textbf{u})^\top(\textbf{X}, \textbf{u})\right)^{1/2}.
\]
Substituting \( \textbf{u} = \sigma^{-1}(\textbf{y} - \textbf{X}\boldsymbol\beta) \), this becomes
\[
J(\textbf{y}, \boldsymbol{\theta}) = n^{-1}\sigma^{-1}\det\left((\textbf{X}, \textbf{y} - \textbf{X}\boldsymbol\beta)^\top(\textbf{X}, \textbf{y} - \textbf{X}\boldsymbol\beta)\right)^{1/2}.
\]
As a consequence of the Cauchy--Binet formula, the determinant of 
\(
(\textbf{X},\textbf{y}-\textbf{X}\boldsymbol\beta)^\top
(\textbf{X},\textbf{y}-\textbf{X}\boldsymbol\beta)
\)
is invariant in \( \boldsymbol\beta \). To see this, write
\[
\textbf{y}-\textbf{X}\beta
=
\textbf{y}-\textbf{X}\hat{\boldsymbol\beta}
+
\textbf{X}(\hat{\boldsymbol\beta}-\boldsymbol\beta),
\qquad
\hat{\boldsymbol\beta}=(\textbf{X}^\top\textbf{X})^{-1}\textbf{X}^\top\textbf{y}.
\]
Since changing \( \boldsymbol\beta \) only modifies the last column by a linear combination of the
columns of \( \textbf{X} \),  the determinant is unchanged.

To simplify computation, we can therefore evaluate the Jacobian at 
\( \boldsymbol\beta=\hat{\boldsymbol\beta}\). At this value, \(\textbf{r}=\textbf{y}-\textbf{X}\hat{\boldsymbol\beta}\) is orthogonal to the columns of 
\(\textbf{X}\), giving
\[
(\textbf{X},\textbf{r})^\top(\textbf{X},\textbf{r})
=
\begin{pmatrix}
\textbf{X}^\top\textbf{X} & 0 \\
0 & \textbf{r}^\top\textbf{r}
\end{pmatrix}.
\]
It follows that
\[
J(\textbf{y}, \boldsymbol{\theta})
=
\sigma^{-1}
\det(\textbf{X}^\top\textbf{X})^{1/2}
\cdot
\text{RSS}^{1/2},
\]
where \( \text{RSS}=\|\mathbf r\|^2 \).

Thus, the density of the GFD is
\[
r(\boldsymbol\beta,\sigma | \textbf{y}) \propto \sigma^{-n-1}p_{\mathbf U}(\sigma^{-1}(\textbf{y}-\textbf{X}\boldsymbol\beta)),
\]
which coincides with the Bayesian solution using a Jeffreys prior,  known to satisfy the Bernstein--von Mises theorem \citep{hannig.2016.tandf}.
\end{example}

\section{A Bernstein--von Mises theorem for generalized fiducial distributions}\label{s:BvMtheorem}
\subsection{Motivation}\label{ss:motivation}
Recall from the previous section that the generalized fiducial density takes a form similar to the Bayesian posterior: \[ r(\boldsymbol{\theta}|\mathbf{y}) \propto p_{\boldsymbol{\theta}}(\mathbf y)J(\mathbf{y,\boldsymbol{\theta}}).\] Here, $J(\mathbf{y,\boldsymbol{\theta}})$ is not pre-specified in any sense, but is determined only by the data model and the distance function in the GFD. It is the volume correction factor resulting from the transformation of the data model to a model on the parameter space, and acts as a weight function on the likelihood, similar to a prior density. In the Bayesian case, a prior density that is continuous and positive is roughly constant in a neighborhood of the true parameter, and so effectively "cancels" from the expression for the posterior density. Under regularity conditions of the model, the likelihood that remains in the expression behaves asymptotically like that of a normal. We discuss such regularity conditions in the next section. 

The same heuristic argument for the asymptotic behavior of a GFD applies if the weight function converges (in $n$) to a function that is, like a well-behaved prior, roughly constant around the truth. For i.i.d. data, this convergence is typically guaranteed. When the \(l_{\infty}\) norm is used in \eqref{eq:pseudoInverse}, the Jacobian function is a U-statistic and it can be shown by applying \citet{yeo2001uniform} to converge a.s. to \( \pi(\boldsymbol{\theta})  = E_{\theta}(J(\mathbf{y,\boldsymbol{\theta}})
)\) uniformly in $\boldsymbol{\theta}$ on compact sets (see \citet{hannig.2009.statistica} for further discussion). Under the \(l_{2}\) norm, \(J(\mathbf{y,\boldsymbol{\theta}})\) is a function of the mean of positive semi-definite matrices. If the DGA is continuous in $\theta$ and the entries of \(J(\mathbf{y,\boldsymbol{\theta}})\) have finite expectations, then \(J(\mathbf{y,\boldsymbol{\theta}})\) will converge to  its expectation. For non-i.i.d. data, this limit exists when a law of large numbers result holds. In this sense, the GFD can be viewed analogously to an empirical Bayes construction, where \(J(\mathbf{y,\boldsymbol{\theta}})\) is an estimator of $\pi(\boldsymbol{\theta})$. As discussed in \citet{hannig.2009.statistica}, the GFD will always be proper, which may or may not be the case for a Bayesian posterior with a (potentially) improper empirical prior.

\subsection{Problem setup and regularity conditions}\label{ss:setup}

In this section, we introduce the regularity conditions on the statistical model needed for local asymptotic normality (LAN). These conditions stipulate that models are smooth enough to permit a quadratic expansion of the log-likelihood ratio, which guarantees that, under a local parameterization, the likelihood ratio process converges asymptotically to that of a Gaussian experiment. Although LAN is a more general theoretical framework, the regularity conditions required are essentially the same as those needed to establish asymptotic normality and efficiency of maximum likelihood estimators in i.i.d. models.

Consider a statistical model $\{P_{\boldsymbol{\theta}}: \boldsymbol{\theta} \in \boldsymbol{\Theta}\}$ defined on a measurable space $(\mathcal{Y}, \mathcal{A})$, where $\boldsymbol{\Theta}\subseteq \mathbb{R}^d$ is open. Suppose we observe a sample $(y_1, \dots, y_n)$ from $P_{\boldsymbol{\theta}}$,  having density $p_{\boldsymbol{\theta}}$ with respect to some dominating measure. The full sample of size $n$ can be viewed as a single observation from the joint model $\{ P_{\boldsymbol{\theta}}^{(n)} : \boldsymbol{\theta} \in \boldsymbol{\Theta}\}.$  In the i.i.d. case, $P_{\boldsymbol{\theta}}^{(n)}$ is a  product measure and we use the representative notation $P_{\boldsymbol{\theta}}^n$ for this case. According to \citet{lecam.1960.uc} (and later \citet{van2000asymptotic}), a sequence of models $\{P^{(n)}_{\boldsymbol{\theta}}: \boldsymbol{\theta} \in \boldsymbol{\Theta}\}$ is said to be \textit{locally asymptotically normal} (LAN) if there exist matrices $\textbf{r}_n$, invertible matrix $I_{\boldsymbol{\theta}}$, and random vectors $\Delta_{n,\boldsymbol{\theta}}$ converging weakly to $N(0,I_{\boldsymbol{\theta}})$ such that, for every converging sequence $\textbf{h}_n \rightarrow \textbf{h}$, 
\[
\ln \frac{dP^{(n)}_{\boldsymbol{\theta}+\textbf{r}_n^{-1}\textbf{h}_n}}{dP^{(n)}_{\boldsymbol{\theta}}} = \textbf{h}^\top\Delta_{n,\boldsymbol{\theta}} - \frac{1}{2}\textbf{h}^\top I_{\boldsymbol{\theta}}\textbf{h} + o_{P^{(n)}_{\boldsymbol{\theta}}}(1), \quad n \to \infty.
\]
The $\Delta_{n,\boldsymbol{\theta}}$ are, in the asymptotic sense, "locally sufficient" statistics. Note that LAN is a very general property that does not require observations to be independent, nor identically distributed. Our main result, given in the next section, is stated under the assumption that the sequence of models $P_{\boldsymbol{\theta}}^{(n)}$ is LAN with norming matrices $\textbf{r}_n = \sqrt{n}I_d.$ In Example~\ref{ex:ar2}, we discuss a class of smooth, weakly dependent models that are LAN with  $\Delta_{n,\boldsymbol{\theta}} = \dfrac{1}{\sqrt{n}}\dot{\ell}_{\boldsymbol{\theta}}(y_{1:n})$, where \(\dot{\ell}_{\boldsymbol{\theta}}(y_{1:n})\) is the full-sample score function.

In the i.i.d. case, LAN is guaranteed by the single-observation differentiability in quadratic mean (DQM) property. This property is much weaker than classical regularity conditions that require twice continuous differentiability of the log-likelihood with respect to $\boldsymbol{\theta}$ and a dominated second derivative, as discussed in Section~\ref{ss:bvm}. To give a precise definition, the family of models $\{P_{\boldsymbol{\theta}}: \boldsymbol{\theta} \in \boldsymbol{\Theta}\}$ is \textit{differentiable in quadratic mean} at $\boldsymbol{\theta} \in \boldsymbol{\Theta}$ if there exists a function $\dot{\ell}_{\boldsymbol{\theta}} : \mathcal{Y} \rightarrow \mathbb{R}^d$ such that 
\[
\int \left( \sqrt{p_{\boldsymbol{\theta}+\textbf{h}}} - \sqrt{p_{\boldsymbol{\theta}}} - \frac{1}{2}\textbf{h}^{\top}\dot{\ell}_{\boldsymbol{\theta}}\sqrt{p_{\boldsymbol{\theta}}}\right)^2 d\mu = o(\norm{\textbf{h}}^2)
\]
as $\textbf{h}\rightarrow \textbf{0}$. Here, $\dot{\ell}_{\boldsymbol{\theta}}$ denotes a generalized derivative (in the quadratic-mean sense), which coincides with the ordinary derivative when it exists. An implication of the DQM condition is that
$\dot{\ell}_{\boldsymbol{\theta}}$ can be viewed as a generalized score function and a generalization of the Fisher information is naturally defined
as its variance, $I_{\boldsymbol{\theta}} = E_{\boldsymbol{\theta}}\dot{\ell}_{\boldsymbol{\theta}}\dot{\ell}_{\boldsymbol{\theta}}^T$. Expanding the log-likelihood ratio using the DQM condition 
produces the LAN expression, with linear coefficient 
$\Delta_{n,\boldsymbol{\theta}} = \dfrac{1}{\sqrt{n}}\sum\limits_{i=1}^{n}\dot{\ell}_{\boldsymbol{\theta}}(y_i)$ 
and quadratic coefficient
$-\frac{1}{2}E_{\boldsymbol{\theta}}\dot{\ell}_{\boldsymbol{\theta}}\dot{\ell}_{\boldsymbol{\theta}}^T$. The important insight is that while the DQM condition requires only a type of one-times differentiability, it is enough to yield the LAN expansion directly, without appeal to the second Bartlett identity needed in the classical setting. In Section~\ref{s:numericalex}, we explore an example and provide numerical illustrations for a distribution where ordinary differentiability fails and the (generalized) Fisher information and LAN property are obtained through DQM.

\subsection{Main result}\label{ss:mainresult}
Suppose we observe data \( \mathbf{Y}= (Y_1, \dots, Y_n)\) generated from a DGA $G(\mathbf{U},\boldsymbol{\theta}_0)$ for the model $P_{\boldsymbol{\theta}_0}$,  where $\boldsymbol{\theta}_0$ is the "true" parameter, which we assume to be a fixed, non-random value in the interior of $\boldsymbol{\Theta}$. Let $\mathbf{y}= (y_1, \dots, y_n)$ denote \textit{any} potential realized sample. Then, to study frequentist asymptotics, we imagine repeatedly sampling \(\mathbf{Y}\) from this DGA and computing the generalized fiducial distribution \(R_{\bar{\boldsymbol{\Theta}}|\mathbf{Y}=\mathbf{y}}\) 
of \(\boldsymbol{\theta}\) for each realized dataset $\mathbf{y}$. Here, just as in the Bayesian Bernstein--von Mises theorem, \(\boldsymbol{\theta}_0\) is a fixed parameter that generated the data; \(R_{\bar{\boldsymbol{\Theta}}|\mathbf{Y}}\) is a random measure, where the randomness is due to the sampling variability of $\mathbf{Y} \sim P_{\boldsymbol{\theta}_0}$. In what follows, probability 
statements and limits are therefore understood with respect to 
$P_{\boldsymbol{\theta}_0}$.  We also recall that the 
generalized fiducial density is given in Eq.~\eqref{eq:GFdensity}, which will be referenced throughout.

\begin{assumption}[Local asymptotic normality]\label{as:dqm}
The sequence of models $\{P^{(n)}_{\boldsymbol{\theta}_0}: \boldsymbol{\theta} \in \boldsymbol{\Theta}\}$ is  locally asymptotically normal around $\boldsymbol{\theta}_0$ with norming matrices $\textbf{r}_n = \sqrt{n}I_d$.
\end{assumption}

\begin{assumption}[Convergence of Jacobian function in GF density]\label{as:limit}There exists a sequence $D_n > 0$ such that $D_n \sqrt{n} \to \infty$, and a $\sigma$-finite measure absolutely continuous with respect to the Lebesgue measure on $\norm{\boldsymbol{\theta}-\boldsymbol{\theta}_0} \leq D_n$ with density $\pi(\boldsymbol{\theta})$, such that as \( n \to \infty \),
  \[\sup\limits_{\boldsymbol{\theta}: \norm{\boldsymbol{\theta}-\boldsymbol{\theta}_0} \leq D_n}\frac{\lvert J(\textbf{Y},\boldsymbol{\theta}) - \pi(\boldsymbol{\theta}) \rvert}{J(\textbf{Y},\boldsymbol{\theta})} \stackrel{P_{\boldsymbol{\theta}_0}}{\rightarrow} 0.\]
\end{assumption}

\begin{assumption}[Mass of limiting measure]\label{as:limitmass} The limiting density $\pi(\boldsymbol{\theta})$ is continuous and positive at $\boldsymbol{\theta}_0$.
\end{assumption}

\begin{assumption}[Likelihood splitting]\label{as:datasplit}
There exist finite measures $P^1_{\boldsymbol{\theta}}, P^2_{\boldsymbol{\theta}}$ absolutely continuous with respect to the Lebesgue measure having densities  $p^1_{\boldsymbol{\theta}}(\textbf{y}), p^2_{\boldsymbol{\theta}}(\textbf{y})$, respectively, such that
    \begin{equation*}
    p_{\boldsymbol{\theta}}(\textbf{y}) = p^1_{\boldsymbol{\theta}}(\textbf{y})p^2_{\boldsymbol{\theta}}(\textbf{y}). 
    \end{equation*}
Moreover, there exists a finite measure on $\theta$, absolutely continuous with respect to the Lebesgue measure, with density $\gamma(\boldsymbol{\theta})$ 
so that for the same $D_n>0$ as in Assumption~\ref{as:limit}, as \( n \to \infty \),
    \begin{equation*}
    P_{\boldsymbol{\theta}_0}\big(p^2_{\boldsymbol{\theta}}(\textbf{Y})J(\textbf{Y},\boldsymbol{\theta})
    I_{\{\norm{\boldsymbol{\theta}-\boldsymbol{\theta}_0} >D_n\}} {\leq} \gamma(\boldsymbol{\theta})\big) \rightarrow 1.
    \end{equation*}
\end{assumption}

\begin{assumption}[Exponentially consistent tests]\label{as:tests}
  For every $\varepsilon >0$, there exists a sequence of tests $\delta_n$ and constants $c_1,c_2>0$ such that for $\norm{\boldsymbol{\theta}-\boldsymbol{\theta}_0} \geq \varepsilon$, 
    \begin{equation*}
        E_{\boldsymbol{\theta}_0}\delta_n \rightarrow 0 \quad\text{and}\quad E_{\boldsymbol{\theta}}(1-\delta_n) \leq e^{-c_1n(\norm{\boldsymbol{\theta}-\boldsymbol{\theta}_0}^2 \wedge 1)} \quad\text{and}\quad 
        E_{\boldsymbol{\theta}}^{1}(1-\delta_n) \leq e^{-c_2n(\norm{\boldsymbol{\theta}-\boldsymbol{\theta}_0}^2 \wedge 1)}.
    \end{equation*}
\end{assumption}
\begin{theorem}\label{thm:BvM}
Under Assumptions 1--5, the fiducial distribution satisfies
\begin{equation*}
\norm{R_{\sqrt{n}(\bar{\boldsymbol{\Theta}}-\boldsymbol{\theta}_0)|\mathbf{Y}} - N(I_{\boldsymbol{\theta}_0}^{-1}\Delta_{n,\boldsymbol{\theta}_0},I_{\boldsymbol{\theta}_0}^{-1})}_{TV} \overset{P_{\boldsymbol{\theta}_0}^{(n)}}{\rightarrow} 0,
\end{equation*} 
where $R_{\sqrt{n}(\bar{\boldsymbol{\Theta}}-\boldsymbol{\theta}_0)|\mathbf{Y}}$ denotes the generalized fiducial distribution on the locally rescaled parameter space and the distance computed is the total variation distance. 
\end{theorem}

In the statement of the theorem, both $R_{\sqrt{n}(\bar{\boldsymbol{\Theta}}-\boldsymbol{\theta}_0)|\mathbf{Y}}$ and $N(I_{\boldsymbol{\theta}_0}^{-1}\Delta_{n,\boldsymbol{\theta}_0},I_{\boldsymbol{\theta}_0}^{-1})$ are measures that are functions of the random variable $\mathbf{Y}$. For each fixed realization of $\mathbf{Y}$, the total variation distance is computed between these two measures, and the theorem states that this (random) distance, viewed as a function of $\mathbf{Y}$, converges to $0$ in probability as $n \to \infty$.

The proof is presented in Section~2 of the Supplementary Material \citep{supplement}; it adapts the well-known Bernstein--von Mises argument presented in \citet{vaart.1998.cambridge} to the generalized fiducial setting. While we require the additional Assumptions~\ref{as:limit} and \ref{as:datasplit} relative to the Bayesian setting, the Bayesian framework starts from the stronger assumption of a fully specified prior density, which GFI does not initially impose. 

As discussed earlier in Section~\ref{ss:motivation}, the local approximation of the Jacobian function as required in Assumption~\ref{as:limit} is expected to hold for a broad class of models. Assumption~\ref{as:datasplit} requires a decomposition of the likelihood into a component used for inference and a component whose role is to control the tail behavior of the Jacobian. This construction is conceptually similar to the minimal training sample approach of \citet{BergerPerichi2004} and the fractional Bayes factor method of \citet{o1995fractional}. In \citet{BergerPerichi2004}, a minimal subset of data that guarantees a proper posterior is used to update an initial (possibly improper) prior, and the resulting proper posterior is then treated as a prior to be used with the rest of the data. The fractional Bayes approach is analogous, in that a fractional power of the likelihood is used to construct a proper posterior, which is then treated as a prior to be used with the remaining fraction of the likelihood for inference. While the GFD is itself proper, a similar construction is useful here for controlling tails in the asymptotic argument. An important point is that Assumption~\ref{as:datasplit} and Assumption~\ref{as:tests} naturally counterbalance each other, meaning a choice of split that makes one assumption vacuously true would likely make it impossible to satisfy the other.

While Assumption~\ref{as:datasplit} is intentionally flexible, a convenient mechanism for verifying this assumption is to introduce an artificial data sample $\mathbf y^{\sharp}.$ This artificial sample does not play a role in the construction of the fiducial procedure, but is introduced strictly as a device to verify Assumption~\ref{as:datasplit}. For a fixed, finite artificial sample large enough so that the Jacobian function exists, the quantity 
$p_{\boldsymbol{\theta}}(\mathbf y^{\sharp})
J(\mathbf y^{\sharp},\boldsymbol{\theta})$ 
is finite by construction. We then define
\[
p_{\boldsymbol{\theta}}^{2}(\mathbf y)
= p_{\boldsymbol{\theta}}(\mathbf y^{\sharp}),
\qquad
p_{\boldsymbol{\theta}}^{1}(\mathbf y)
= \frac{p_{\boldsymbol{\theta}}(\mathbf y)}
       {p_{\boldsymbol{\theta}}(\mathbf y^{\sharp})},
\]
where $\mathbf y^{\sharp}$ denotes the artificial data sample. Multiplying and dividing the Jacobian by 
$J(\mathbf y^{\sharp},\boldsymbol{\theta})$ yields
\[
p_{\boldsymbol{\theta}}^{2}(\mathbf y)
J(\mathbf y,\boldsymbol{\theta})
=
p_{\boldsymbol{\theta}}(\mathbf y^{\sharp})
J(\mathbf y^{\sharp},\boldsymbol{\theta})
\frac{J(\mathbf y,\boldsymbol{\theta})}
     {J(\mathbf y^{\sharp},\boldsymbol{\theta})}.
\]
It therefore remains to show that the ratio
\(
J(\mathbf y,\boldsymbol{\theta})/
     J(\mathbf y^{\sharp},\boldsymbol{\theta})
\)
can be bounded with high probability under 
$P_{\boldsymbol{\theta}_0}$ using arguments similar to those required to 
verify Assumption~\ref{as:limit}. Moreover, because $p_{\boldsymbol{\theta}}(\mathbf y^{\sharp})$ is fixed and finite, exponentially consistent tests that hold under $P_{\boldsymbol{\theta}}$ readily extend to $P_{\boldsymbol{\theta}}^{1}$. Note that $P_{\boldsymbol{\theta}}^{1}$ need not be a probability measure; Assumption~\ref{as:tests} requires only that $E_{\boldsymbol{\theta}}^{1}(1-\delta_n)$ be exponentially bounded. 

As in the Bayesian framework, in the i.i.d. setting Assumption~\ref{as:tests} may be replaced by the weaker requirement that uniformly consistent estimators exist, since such estimators can be used to construct uniformly exponentially consistent tests (see Lemma 10.6 in \citet{vaart.1998.cambridge}). Beyond the i.i.d. case, \citet{ghosal2007convergence} derived conditions for the existence of exponentially consistent tests in stationary Gaussian time series models, an example of which we study in Example~\ref{ex:ar2}. More recently, \citet{connault.2021.working} defined conditions under which uniformly exponentially consistent tests can be constructed from uniformly consistent estimators for a broader class of weakly dependent models.
 
Before presenting numerical examples, we demonstrate in the next section how our result extends theoretical guarantees of GFDs for a class of models previously studied in \citet{sonderegger.2014.springer} and \citet{sonderegger.dissertation}. 

\subsection{Fiducial theory for free-knot spline models}\label{ss:splineex}
\citet{sonderegger.2014.springer} derived a generalized fiducial solution to the free-knot spline problem and evaluated its performance against the Bayesian approach developed by \citet{dimatteo.2001.oxford}, based on a prior of $\pi(\boldsymbol{\alpha},\mathbf{t},\sigma^2) \propto \sigma^{-2}$, where $\sigma^2$ is a variance term, $\boldsymbol{\alpha}$ denotes the polynomial coefficients, and $\mathbf{t}$ represents the vector of knot points. In that analysis, the coverage rates of fiducial intervals were generally closer to the nominal levels, while the Bayesian approach produced credible intervals that either under-covered compared to the fiducial or were farther from the nominal level when both methods were conservative. \citet{sonderegger.2014.springer} also demonstrated that splines of degree $p\geq4$ satisfy conditions for asymptotic normality (based on more restrictive differentiability criteria than used here) of their multivariate generalized fiducial estimator. Moreover, simulation studies of the same solution given in Chapters 4 and 5 of \citet{sonderegger.dissertation} indicate that the fiducial solution for the free-knot spline problem is asymptotically accurate and demonstrates strong performance for lower degree splines (i.e. degree $p\leq3 $). Importantly, lower degree splines are more realistic in practical applications, so extending the theory to these cases is meaningful. 

It is well-known that exponential family models are differentiable in quadratic mean, and here, we show that splines of degree $p\geq 1$ satisfy the other conditions of Theorem~\ref{thm:BvM}, supporting the numerical results provided in \citet{sonderegger.2014.springer} and \citet{sonderegger.dissertation}. Because extensive simulation studies are already presented in those earlier papers, we do not duplicate them here and instead focus on extending the theory. We provide a detailed derivation of the GFD and verify the remaining conditions of the theorem as an illustrative example; in subsequent applications, analogous arguments are presented more concisely and are deferred to Sections~3 and 4 of the Supplementary Material \citep{supplement}.

For background, a spline model of degree $p$ is characterized by piecewise degree $p$ polynomials with the requirement that the resulting function be, in some sense, smooth at the connection, or \textit{knot}, points. In the \emph{free-knot} setting, the number of knots $\kappa$ is known but their locations $\mathbf{t}=(t_1,\ldots,t_{\kappa})^\top$ are unknown and to be estimated \citep{toms.2003.ecology,sonderegger.2009.frontiers}. Consider an observed data sample $\{(x_i,y_i)\}_{i=1}^n$ generated from
\[
Y_i = g(x_i \mid \boldsymbol{\alpha},\mathbf{t}) + \sigma \varepsilon_i, 
\quad \varepsilon_i \stackrel{\text{iid}}{\sim} N(0,1),
\]
where $x_i \in [a,b]$ are observed design points. We treat the design points as observed realizations from a distribution with positive density on $[a,b]$, and all inference is conditional on the observed $x_i$. The spline has coefficients $\boldsymbol{\alpha}=(\alpha_0,\ldots,\alpha_{p+\kappa})^\top$ and knot vector $\mathbf{t}$ satisfying
\[
t_{k+1} > t_k, \quad |t_{k+1}-t_k| > \delta, \quad 
t_1 \ge a+\delta, \quad t_\kappa \le b-\delta,
\]
where $\delta>0$, and $|\alpha_{p+k}| > \xi > 0$ for identifiability.

The truncated polynomial representation is
\[
g(x_i \mid \boldsymbol{\alpha},\mathbf{t}) 
= \sum_{j=0}^{p}\alpha_jx_i^j 
+ \sum_{k=1}^{\kappa}\alpha_{p+k}(x_i-t_k)_+^p,
\]
where $(z)_+=\max\{z,0\}$. The spline basis is given by
\[
B(x,\mathbf{t})
= \{1,x,\ldots,x^p,(x-t_1)^p_+,\ldots,(x-t_{\kappa})^p_+\},
\]
and the conditional density can then be written as
\[
p_{\boldsymbol{\theta}}(y_i\mid x_i)
= \frac{1}{\sqrt{2\pi\sigma^2}}
\exp\!\Big(-\frac{1}{2\sigma^2}
\big(y_i - B(x_i,\mathbf{t})\boldsymbol{\alpha}\big)^2\Big),
\]
with parameter vector $\boldsymbol{\theta}=(\mathbf{t},\boldsymbol{\alpha},\sigma^2)^\top$.

A suitable data generating algorithm for constructing a GFD is then
\[
Y_i = G(\varepsilon_i,\boldsymbol{\theta}) 
= g(x_i\mid\boldsymbol{\alpha},\mathbf{t}) + \sigma \varepsilon_i.
\]
Given the observed sample $\{y_i\}_{i=1}^n$, inversion yields \[\varepsilon_i = \frac{1}{\sigma}(y_i - g(x_i\mid\boldsymbol{\alpha},\mathbf{t})),\]
and differentiating $G$ with respect to $\boldsymbol{\theta}$ gives
\begin{align*}
\frac{\partial G}{\partial \boldsymbol{\alpha}}
&= \big(1,x_i,\ldots,x_i^p,(x_i-t_1)_+^p,\ldots,(x_i-t_{\kappa})_+^p\big)\\
\frac{\partial G}{\partial \mathbf{t}}
&= p\big(\alpha_{p+1}(x_i-t_1)_+^{p-1},\ldots,\alpha_{p+\kappa}(x_i-t_{\kappa})_+^{p-1}\big), \\
\frac{\partial G}{\partial \sigma}
&= \frac{1}{\sigma}(y_i - g(x_i\mid\boldsymbol{\alpha},\mathbf{t})).
\end{align*}
From these, we form the gradient matrix $\mathbf{A}=[\mathbf{B}_{\boldsymbol{\alpha}},\mathbf{B}_{\mathbf{t}},\mathbf{B}_{\sigma}]$, where
\[
\mathbf{B}_{{\sigma}}=(y_1 - g(x_1\mid\boldsymbol{\theta}),\ldots,y_n - g(x_n\mid\boldsymbol{\theta}))^\top,
\]
and  
\[
\mathbf{B}_{\boldsymbol\alpha} = 
\begin{pmatrix}
1 & x_1 & \dots & x_1^p & (x_1-t_1)_+^p & \dots & (x_1-t_\kappa)_+^p \\
\vdots & \vdots & & \vdots & \vdots & & \vdots \\
1 & x_n & \dots & x_n^p & (x_n-t_1)_+^p & \dots & (x_n-t_\kappa)_+^p
\end{pmatrix}_{n \times (p+\kappa+1)},
\]
\[
\mathbf{B}_{\mathbf{t}} =
\begin{pmatrix}
(x_1-t_1)_+^{p-1} & \dots & (x_1-t_\kappa)_+^{p-1} \\
\vdots & & \vdots \\
(x_n-t_1)_+^{p-1} & \dots & (x_n-t_\kappa)_+^{p-1}
\end{pmatrix}_{n \times \kappa}.
\]
Since $\mathbf{B}_{\sigma}$ is a linear combination of columns of $\mathbf{B}_{\boldsymbol{\alpha}}$, it can be replaced (in the determinant) by $\mathbf{y}= (y_1,\ldots,y_n)^\top$. The Jacobian is therefore
\[
J(\mathbf{y},\boldsymbol{\theta})
= \frac{p^{\kappa}}{\sigma}
\prod_{k=1}^{\kappa}|\alpha_{p+k}|
\det\!\Big[\frac{1}{n}\tilde{\mathbf{A}}^\top\tilde{\mathbf{A}}\Big]^{1/2},
\quad \tilde{\mathbf{A}}=[\mathbf{B}_{\boldsymbol{\alpha}},\mathbf{B}_{\mathbf{t}},\mathbf{y}].
\]

The Fisher information is well-defined without requiring the more general DQM formulation, and is non-singular. The score for $\sigma^2$ is
\[
\frac{\partial}{\partial\sigma^2}\ell_{\boldsymbol{\theta}}(y_i)
= -\frac{1}{2\sigma^2}
+ \frac{1}{2(\sigma^2)^2}\big(y_i - B(x_i,\mathbf{t})^\top\boldsymbol{\alpha}\big)^2,
\]
and direct calculation gives
\[
E_{\boldsymbol{\theta}}\!\left[
\left(\frac{\partial \ell}{\partial\sigma^2}\right)^2
\right]
= \frac{1}{2\sigma^4} > 0.
\]
Moreover, for $\boldsymbol{\eta} = (\boldsymbol{\alpha},\mathbf{t})$, we have
\[
E_{\boldsymbol{\theta}}\!\left[
\frac{\partial \ell}{\partial\eta_j}
\frac{\partial \ell}{\partial\sigma^2}
\right]
= 0,
\]
since $\partial \ell/\partial \eta_j$ is linear in the residual $\varepsilon_i = y_i - B(x_i,\mathbf{t})^\top \boldsymbol{\alpha}$ and $\partial \ell/\partial \sigma^2$ involves only constants and $\varepsilon_i^2$. Hence, the Fisher information decomposes into blocks corresponding to $\boldsymbol{\eta}$ and $\sigma^2$.

For $\boldsymbol{\eta}$,
\[
I_{jk}(\boldsymbol{\eta})
= E_{\boldsymbol{\theta}}\bigg[\frac{\partial \ell}{\partial\eta_j}\frac{\partial \ell}{\partial\eta_k}\bigg]  
   =\frac{1}{\sigma^2}\frac{\partial}{\partial\eta_j}B(x_i,\mathbf{t})^\top\boldsymbol{\alpha}\frac{\partial}{\partial\eta_k}B(x_i,\mathbf{t})^\top\boldsymbol{\alpha},
\]
where
\[
\frac{\partial}{\partial t_j}B(x_i,\mathbf{t})^\top\boldsymbol{\alpha}
= -p(x_i-t_j)^{p-1}_+\alpha_{p+j},
\quad
\frac{\partial}{\partial\alpha_j}B(x_i,\mathbf{t})^\top\boldsymbol{\alpha}
= B_j(x_i,\mathbf{t}).
\]
Let $\mathbf{X}$ denote the $n\times(p+2\kappa+1)$ matrix with columns
\[
\mathbf{X}
= 
\big[B_1(x_i,\mathbf{t}),\ldots,B_{p+\kappa+1}(x_i,\mathbf{t}),
\frac{\partial}{\partial t_1}B(x_i,\mathbf{t})^\top\boldsymbol{\alpha},\ldots,
\frac{\partial}{\partial t_\kappa}B(x_i,\mathbf{t})^\top\boldsymbol{\alpha}\big].
\]
Then the Fisher information matrix is
\[
I_{\boldsymbol{\theta}} =
\begin{pmatrix}
I(\boldsymbol{\eta}) & 0 \\
0 & 1/(2\sigma^4)
\end{pmatrix}, 
\quad
I(\boldsymbol{\eta}) = \frac{1}{\sigma^2}\mathbf{X}^\top \mathbf{X},
\]
which is non-singular if the columns of $\mathbf{X}$ are linearly independent, i.e., when $t_j \neq t_k$ and $\alpha_{p+j}\neq 0$ for all $j$.

Using notation and expressions derived above, we can write the score for a single observation $y_i$ as
\[
\dot{\ell}_{\boldsymbol{\theta}}(y_i)
=
\begin{pmatrix}
\frac{1}{\sigma^2} \frac{\partial}{\partial \boldsymbol{\eta}} B(x_i,\mathbf{t})^\top \boldsymbol{\alpha} \,(y_i - B(x_i,\mathbf{t})^\top \boldsymbol{\alpha}) \\
-\frac{1}{2\sigma^2} + \frac{1}{2\sigma^4} (y_i - B(x_i,\mathbf{t})^\top \boldsymbol{\alpha})^2
\end{pmatrix},
\quad \boldsymbol{\eta} = (\boldsymbol{\alpha}, \mathbf{t}).
\]

Let $\mathbf{r} = (y_1 - B(x_1,\mathbf{t})^\top \boldsymbol{\alpha}, \dots, y_n - B(x_n,\mathbf{t})^\top \boldsymbol{\alpha})^\top$ denote the residuals. Then,
\[
I_{\boldsymbol{\theta}}^{-1}\Delta_{n,\boldsymbol{\theta}}
= \frac{1}{\sqrt{n}} \sum_{i=1}^{n} I_{\boldsymbol{\theta}}^{-1} \dot{\ell}_{\boldsymbol{\theta}}(y_i)
=
\begin{pmatrix}
\frac{1}{\sqrt{n}} (\mathbf{X}^\top \mathbf{X})^{-1} \mathbf{X}^\top \mathbf{r} \\[0.5em]
\frac{1}{\sqrt{n}} \big( \mathbf{r}^\top\mathbf{r} - n\sigma^2 \big)
\end{pmatrix}
\]

and we can verify the remaining conditions of Theorem~\ref{thm:BvM} to show that
\[
\norm{R_{\sqrt{n}(\bar{\boldsymbol{\Theta}}-\boldsymbol{\theta}_0)\mid \mathbf{Y}} 
- N\big(I_{\boldsymbol{\theta}_0}^{-1}\Delta_{n,\boldsymbol{\theta}_0}, I_{\boldsymbol{\theta}_0}^{-1}\big)}_{TV} 
\;\;\overset{P_{\boldsymbol{\theta}_0}^{n}}{\longrightarrow} 0.
\]

\textit{Assumption \ref{as:limit} (Convergence of Jacobian function)} and \textit{Assumption \ref{as:limitmass} (Mass of limiting measure).}  
Note that
\[
\frac{1}{n} \tilde{\mathbf{A}}^\top \tilde{\mathbf{A}} =
\frac{1}{n} 
\begin{pmatrix}
\mathbf{B}_{\boldsymbol\alpha}^\top \mathbf{B}_{\boldsymbol\alpha} & \mathbf{B}_{\boldsymbol\alpha}^\top \mathbf{B}_{\mathbf t} & \mathbf{B}_{\boldsymbol\alpha}^\top \mathbf{y} \\
\mathbf{B}_{\mathbf t}^\top \mathbf{B}_{\boldsymbol\alpha} & \mathbf{B}_{\mathbf t}^\top \mathbf{B}_{\mathbf t} & \mathbf{B}_{\mathbf t}^\top \mathbf{y} \\
\mathbf{y}^\top \mathbf{B}_{\boldsymbol\alpha} & \mathbf{y}^\top \mathbf{B}_{\mathbf t} & \mathbf{y}^\top \mathbf{y}
\end{pmatrix},
\] and each element involves sums of terms of the form 
\(x_i^j(x_i - t_k)_+^{\ell}y_i^m\). 
Since \((x_i,y_i)\) are i.i.d.\ and the basis functions are bounded on \([a,b]\), the uniform law of large numbers implies the Jacobian function converges as
\[
J(\mathbf{y},\boldsymbol{\theta})
\;\xrightarrow[]{a.s.}\;
\pi(\boldsymbol{\theta})
= \frac{p^{\kappa}}{\sigma}
\Bigg(\prod_{k=1}^{\kappa}\!|\alpha_{p+k}|\!\Bigg)
\det\!\Bigg(E_{\boldsymbol \theta_0}\!\left[\tfrac{1}{n}\tilde{\mathbf{A}}^\top\tilde{\mathbf{A}}\right]\!\Bigg)^{1/2},
\]
where the expectation is taken elementwise. 
Because all basis terms are bounded and each \(\alpha_{p+k}\) is bounded away from zero, 
\(\pi(\boldsymbol{\theta})\) is finite and positive, satisfying 
Assumptions~\ref{as:limit} and~\ref{as:limitmass}.

\medskip
\noindent
\textit{Assumption \ref{as:datasplit} (Likelihood splitting).}  
For this, we construct a minimal training sample by generating $m = p+\kappa+2$ artificial data points $\{x_i^\sharp\}$, with at least two points between consecutive knots and $|x_i^\sharp - x_j^\sharp| > \delta/2$ for $i\neq j$, and set $y_i^\sharp \in [-q,q]$. Let $\tilde{\mathbf{A}}^\sharp$ be the corresponding gradient matrix in the fiducial solution and $J(\mathbf{y}^\sharp,\boldsymbol{\theta})$ the Jacobian term; by construction and $\alpha_{p+k}>\xi>0$, this Jacobian is bounded away from zero.  

Define
\[
p^2_{\boldsymbol{\theta}}(\mathbf{y}) = p_{\boldsymbol{\theta}}(\mathbf{y}^\sharp) = \prod_{i=1}^m \frac{1}{\sqrt{2\pi\sigma^2}}\exp\Big[-\frac{1}{2\sigma^2}(y_i^\sharp-B(x_i^\sharp,\mathbf{t})\boldsymbol{\alpha})^2\Big], 
\quad
p^1_{\boldsymbol{\theta}}(\mathbf{y}) = \frac{p_{\boldsymbol{\theta}}(\mathbf{y})}{p_{\boldsymbol{\theta}}(\mathbf{y}^\sharp)}.
\]

\noindent Then
\[
p^2_{\boldsymbol{\theta}}(\mathbf{y}) J(\mathbf{y},\boldsymbol{\theta})
= p_{\boldsymbol{\theta}}(\mathbf{y}^\sharp) J(\mathbf{y}^\sharp,\boldsymbol{\theta}) 
\frac{\det\big[\frac{1}{n}\tilde{\mathbf{A}}^\top \tilde{\mathbf{A}}\big]}{\det\big[\frac{1}{m}\tilde{\mathbf{A}}^{\sharp^\top} \tilde{\mathbf{A}}^\sharp\big]}.
\]
\noindent
By construction of the artificial points, for $\mathbf{t}$ satisfying the knot constraints,
\[
\underset{\mathbf{t}}{\inf} \det\Big[\frac{1}{m} \tilde{\mathbf{A}}^{\sharp^\top} \tilde{\mathbf{A}}^\sharp \Big] > 0,
\quad
\underset{\mathbf{t}}{\sup} \det\Big[\frac{1}{m} \tilde{\mathbf{A}}^{\sharp^\top} \tilde{\mathbf{A}}^\sharp \Big] < C_1,
\]
since the $x_i^\sharp$ are separated and $\alpha_{p+k}>\xi>0$.  

Moreover, this term in the Jacobian function for the observed data satisfies a uniform law of large numbers:
\[
\frac{1}{n} \tilde{\mathbf{A}}^\top \tilde{\mathbf{A}} \;\longrightarrow\; E_{\boldsymbol \theta_0}\Big[\frac{1}{n} \tilde{\mathbf{A}}^\top \tilde{\mathbf{A}}\Big], \quad \text{uniformly over all $\mathbf{t}$ satisfying the knot constraints},
\]
so for any $\nu>0$,
\[
P_{\boldsymbol \theta_0}\Bigg(\sup_{\mathbf{t}\in T_\delta} \Big\| \frac{1}{n}\tilde{\mathbf{A}}^\top \tilde{\mathbf{A}} - E_{\boldsymbol \theta_0}\Big[\frac{1}{n} \tilde{\mathbf{A}}^\top \tilde{\mathbf{A}}\Big] \Big\|_F < \nu \Bigg) \longrightarrow 1,
\]
where $\|\cdot\|_F$ is the Frobenius norm.  

Combining these facts, there exists a constant $C_2>0$ such that
\[
P_{\boldsymbol \theta_0}\Big(p^2_{\boldsymbol{\theta}}(\mathbf{y}) J(\mathbf{y},\boldsymbol{\theta}) < C_2 \, p_{\boldsymbol{\theta}}(\mathbf{y}^\sharp) J(\mathbf{y}^\sharp,\boldsymbol{\theta})\Big) \longrightarrow 1.
\]

Finally, since $\int p_{\boldsymbol{\theta}}(\mathbf{y}^\sharp) J(\mathbf{y}^\sharp,\boldsymbol{\theta}) \, d\boldsymbol{\theta} < \infty$, we can set
\(
\gamma(\boldsymbol{\theta}) = C_2 \, p_{\boldsymbol{\theta}}(\mathbf{y}^\sharp) J(\mathbf{y}^\sharp,\boldsymbol{\theta}).
\)

\textit{Assumption \ref{as:tests} (Exponentially consistent tests).}  
\noindent
For the original model $P_{\boldsymbol{\theta}}$, the exponential bound follows directly from Lemma~10.6 in \citet{vaart.1998.cambridge}, so Assumption~\ref{as:tests} holds.  

\noindent
It remains to verify the exponential bound for $P_{\boldsymbol{\theta}}^1$. To this end we decompose $p_{\boldsymbol{\theta}}^1$, separating the first $m < M$ data points:
\begin{align*}
p^1_{\boldsymbol{\theta}}(\mathbf{y})
&= \left(\frac{1}{\sqrt{2\pi\sigma^2}}\right)^{M-m}\exp\Big\{\frac{1}{2\sigma^2}\Big(\sum_{i=1}^m (y_i^\sharp - B(x_i^\sharp,\mathbf{t})^\top\boldsymbol{\alpha})^2 - \sum_{i=1}^M (y_i - B(x_i,\mathbf{t})^\top\boldsymbol{\alpha})^2\Big)\Big\}\\
&\quad\times \prod_{i=M+1}^n \frac{1}{\sqrt{2\pi\sigma^2}} \exp\Big\{-\frac{1}{2\sigma^2}(y_i - B(x_i,\mathbf{t})^\top\boldsymbol{\alpha})^2\Big\}.
\end{align*}

\noindent
The first term is bounded because the difference of quadratic forms is negative definite for a finite $M$ and $\boldsymbol{t}$ in a compact set. The second term is a subset of the original model $P_{\boldsymbol{\theta}}$ with exponential bound $e^{-b(n-M)\norm{\boldsymbol{\theta}-\boldsymbol{\theta}_0}^2} \leq e^{-cn\norm{\boldsymbol{\theta}-\boldsymbol{\theta}_0}^2} $ for $c \geq {b(n-M)}/{n}.$

\section{Examples and simulations}\label{s:numericalex}
We conclude the paper with two interesting examples to further illustrate the generality of our result in models where classical regularity conditions fail or i.i.d. assumptions are not met.
\begin{example}
As an example that fails classical regularity conditions, we consider the triangular distribution. The probability density function for this  distribution is given by

\vspace{-0.5cm}

\begin{equation}\label{eq:triangledensity}
p_{{\theta}}(y) = \dfrac{2y}{\theta}I_{(0,\theta]}(y) + \dfrac{2-2y}{1-\theta}I_{(\theta,1)}(y),\quad \theta\in(0,1).
\end{equation}
The special case $\theta = 1/2$ refers to the \textit{symmetric} triangular distribution, which corresponds to the distribution of the mean of two i.i.d. uniform random variables. While the map $\theta \mapsto p_{{\theta}}(y)$ is continuous, it is not differentiable at $\theta = y$. Indeed, the one-sided derivatives disagree at that point, which violates the pointwise differentiability (in $\theta$) required by classical regularity conditions. However, since differentiability fails only at a single point, it is differentiable in quadratic mean, which requires a derivative only in mean-square rather than at every point (see Section~3 of the Supplementary Material for verification \citep{supplement}). As discussed in Section~\ref{ss:setup}, DQM directly defines the (generalized) Fisher information via the generalized score function, $I_{{\theta}} = E_{{\theta}}\dot{\ell}_{{\theta}}\dot{\ell}_{{\theta}}^T$, without requiring pointwise differentiability of the log-likelihood. The (generalized) Fisher information for the triangular distribution can then be computed using this definition, with (generalized) score function 
\[\dot{\ell}_{\theta}(y_i) = -1/\theta \, I_{(0,\theta]}(y_i) + 1/(1-\theta) \, I_{(\theta,1)}(y_i),\] which gives,
\[
I_{\theta}
= \int_0^{\theta}\!\!\left(-\tfrac{1}{\theta}\right)^2\tfrac{2y}{\theta}dy
+ \int_{\theta}^{1}\!\!\left(\tfrac{1}{1-\theta}\right)^2\tfrac{2-2y}{1-\theta}dy
= \tfrac{1}{\theta} + \tfrac{1}{1-\theta}
= \big(\theta(1-\theta)\big)^{-1}, \quad \theta \in (0,1).
\]

From this, we can formulate the Jeffreys prior  \citep{jeffreys.1998.theory}, which is proportional to the positive square root of the determinant of the Fisher information matrix. \citet{berger.2009.annals} conjectured that the Jeffreys prior for the triangular distribution should be the $\text{Beta}(1/2, 1/2)$ distribution based on numerical evaluation, but could not derive it analytically, citing lack of differentiability. The derivation of  the Fisher information based on DQM shows that the Jeffreys prior is indeed proportional to $I^{1/2}_{\theta} = (\theta(1-\theta))^{-1/2}$, confirming the conjecture in \citet{berger.2009.annals}. 

We obtain the generalized fiducial solution by constructing a suitable DGA from the inverse of the cumulative distribution function (CDF) of the triangular distribution:
\[
Y_i = G(U_i, \theta)=  \sqrt{U_i\theta}\, I_{(0,\theta]}(U_i)
+ \left[1 - \sqrt{(1-U_i)(1-\theta)}\,\right] I_{(\theta,1)}(U_i),
\quad i=1,\ldots,n,
\]
where $\mathbf{U}=(U_1,\ldots,U_n)$ are i.i.d.\ $\text{Uniform}(0,1)$ variables.

The corresponding Jacobian function for this DGA is
\[
J(\mathbf{y},\theta)
= \left[
\frac{1}{n}\sum_{i=1}^{n}
\left\{
\left(\frac{y_i}{2\theta}\right)^2 I_{(0,\theta]}(y_i)
+ \left(\frac{1-y_i}{2(1-\theta)}\right)^2 I_{(\theta,1)}(y_i)
\right\}
\right]^{1/2}.
\]

The remaining conditions for Theorem~\ref{thm:BvM} are verified in Section~3 of the Supplementary Material \citep{supplement}, showing that for a true parameter $\theta_0  \in (0,1)$, 
\[
\norm{R_{\sqrt{n}(\bar{\Theta}-\theta_0)\mid \mathbf{Y}} 
- N\big(I_{\theta_0}^{-1}\Delta_{n,\theta_0}, I_{\theta_0}^{-1}\big)}_{TV} 
\;\;\overset{P_{\boldsymbol{\theta}_0}^{n}}{\longrightarrow} 0,
\]

where the inverse Fisher information is
\(
I_{\theta_0}^{-1} = \theta_0(1-\theta_0)
\)
and the mean is
\[
I_{\theta_0}^{-1}\Delta_{n,\theta_0} = \frac{1}{\sqrt{n}} \sum_{i=1}^{n} I_{\theta_0}^{-1} \dot{\ell}_{\theta_0}(y_i) 
= \frac{n\theta_0 - \sum_{i=1}^{n} I_{(0,\theta_0]}(y_i)}{\sqrt{n}}.
\]

Asymptotically normal behavior of both the GFD and Bayesian posteriors under two non-informative priors is illustrated in Figure~\ref{fig:tri_densities}. Each panel displays density curves for the GFD (red), the Bayesian posterior with Jeffreys prior (blue), and the Bayesian posterior with a flat prior (black) rescaled to the local scale, all based on samples of the size indicated for that panel generated with a true parameter value of $\theta_0 = 0.3$. The density curve from the corresponding normal approximation [$N\big(I_{\theta_0}^{-1}\Delta_{n,\theta_0}, I_{\theta_0}^{-1}\big)$] is overlaid in each panel as a black dashed curve. As the sample size increases (left to right), the GFD and Bayesian posterior density curves align with the normal. 

\vspace{-0.5cm}

\begin{figure}[H]
	\centering
\includegraphics[width=0.85\linewidth]{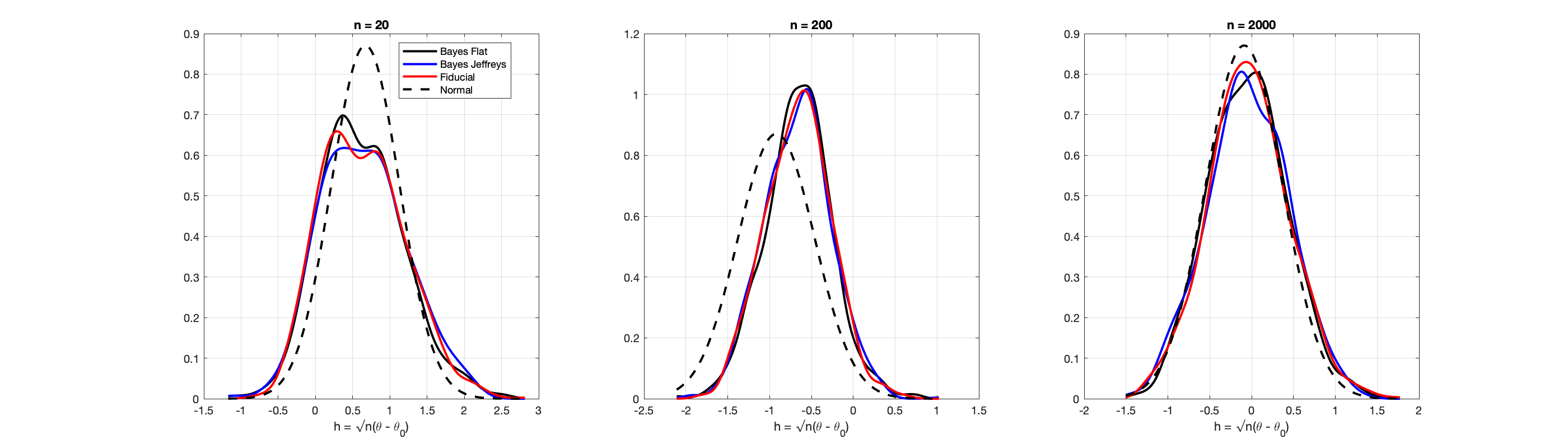}
	\caption{Density curves for the GFD (red) and Bayesian posteriors with Jeffreys (blue) and flat (black) priors, based on samples of size $n= 20, 200, 2,000$ (left to right) for a true parameter value of $\theta_0 = 0.3$. Densities are shown on the local scale, with corresponding normal approximations overlaid (black dashed). As $n$ increases, all distributions converge to the normal limit.}
    \label{fig:tri_densities}
\end{figure}

\vspace{-1cm}

In the single parameter case, a DGA that is one-to-one in \(\theta\) for each \(u\) will lead to a GFD that has exact coverage of \(\theta_0\) at \(n=1\). However, if  there is no \(\theta\) in the parameter space solving \(y = g(u,\theta),\)
(i.e., the inverse-image would produce values outside the parameter space), then those \((u, y)\) pairs are ignored and contribute no mass to the fiducial measure. This is the case for the present example, where the DGA for the triangular values has inverse map:

\[ G^{-1}(u, y) = \begin{cases} 
      \varnothing & u< y^2 \\
      \frac{y^2}{u} & y^2\leq u< y \\
      1 - \frac{(1-y)^2}{1-u} & y\leq u \leq 1-(1-y)^2 \\
      \varnothing & 1-(1-y)^2 < u,
   \end{cases}
\]
where $\varnothing$ denotes the empty set. 

\begin{figure}[H]
    \centering
    \includegraphics[width=0.95\linewidth]{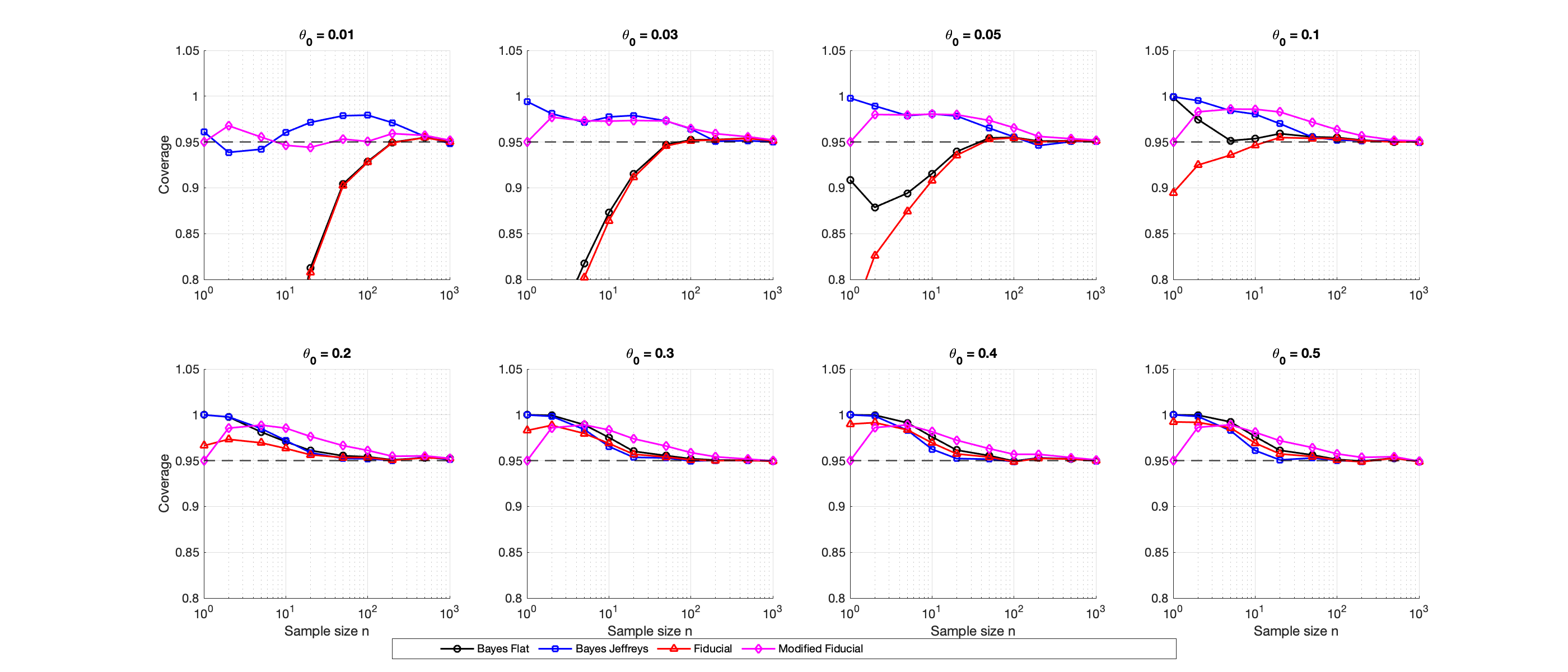}
    \caption{Coverage probabilities across sample sizes for various \(\theta_0\) (blue: Jeffreys Bayesian; magenta: modified GFD; black: flat Bayesian; red: GFD), illustrating improved boundary performance of Jeffreys Bayesian and modified GFD methods, with all approaches converging for large samples.}
    \label{fig:triangle_coverage}
\end{figure}

\vspace{-0.5cm}

The construction of the standard GFD removes values of $u$ for which $y = G(\cdot,\theta)$ has no solution and re-normalizes the probabilities, which leads to under-coverage when the truth \((\theta_0)\) is at or near the boundary. This phenomenon is observed in general in cases where the parameter space is in some way constrained; in those cases, \citet{hannig.2016.tandf} (Section 2.1) recommends an exception to the standard GFD construction and instead suggests extending the parameter space, performing the inversion in the extended space, and then projecting to the boundary of the constrained parameter space. 

This under-coverage issue is not unique to GFI, but can arise in any likelihood-based inferential procedure when the true parameter lies at or near the boundary of a constrained space. For example, in the triangular distribution, a Bayesian posterior proportional to the likelihood (i.e. using a flat prior) suffers from the same redistribution to the interior as the GFD, causing under-coverage near the boundary. This under-covering issue at the boundary is illustrated in the top row of Figure~\ref{fig:triangle_coverage}. Each panel reports the empirical coverage, based on $10,000$ replicates, of nominal $95\%$ two-sided generalized fiducial intervals and Bayesian posterior credible intervals across a range of sample sizes $n$, for a fixed true parameter value $\theta_0$. The panels in the top row of that figure correspond to parameter values near the boundary. The red and black lines correspond to the GFD approach and the Bayesian approach using a flat (uniform) prior, respectively, and are shown to under-cover at most sample sizes.  For very large values of $n$, accurate coverage is achieved by both approaches, as expected from the corresponding Bernstein--von Mises results.

An appropriate frequentist correction to this issue is to treat evidence for values outside the parameter space as evidence that the true parameter lies at the boundary, and assign probability to the boundary points (essentially, the extend-and-project idea suggested in \citet{hannig.2016.tandf}). In the Bayesian case, the Jeffreys prior is guaranteed to be first-order probability matching in the interior, and though not optimal at the boundaries, still puts relatively more weight towards the boundaries than a flat prior. Probability matching GFDs and exact coverage in single parameter models are discussed in \citet{majumder2016higher}. In this example, following the recommendation in \citet{hannig.2016.tandf} in the fiducial case, we obtain the below modification of the inverse image (leading to exact coverage for $n=1$):

\vspace{-0.5cm}

\[ \tilde{G}^{-1}(u, y) = \begin{cases} 
      1 & u< y^2 \\
      \frac{y^2}{u} & y^2\leq u< y \\
      1 - \frac{(1-y)^2}{1-u} & y\leq u \leq 1-(1-y)^2 \\
      0 & 1-(1-y)^2 < u.
   \end{cases}
\]

\vspace{-0.5cm}

For $n>1$, a boundary correction is applied whenever any observation $y_i$ would place $\theta$ outside the interior. Thus, the smallest $y_i^2$ for the upper boundary and smallest $(1-y_i)^2$ for the lower boundary determine a threshold for placing $\theta$ at 1 or 0. Because the $u_i$ are uniform realizations, these values directly determine the probability mass assigned to the boundaries. The remaining probability over the interior is assigned to the original GFD. Formally, this gives

\vspace{-0.5cm}

\[
r^{\mathrm{mod}}(\theta \, |\,\mathbf{y}) = 
\min_i (1-y_i)^2 \,\delta_{\{0\}} +
\min_i y_i^2 \,\delta_{\{1\}} +
\Big(1 - \min_i y_i^2 - \min_i (1-y_i)^2\Big) r(\theta\,|\, \mathbf{y}),
\]
where $r(\theta\,|\, \mathbf{y})$ is the original (unmodified) GFD. 

Asymptotic normality also holds for the modified GFD since
\[
\norm{R^{\text{mod}}_{\boldsymbol{\theta} | \mathbf{y}} - R_{\boldsymbol{\theta} | \mathbf{y}}}_{TV}
\leq \min_i y_i^2 + \min_i (1-y_i)^2,
\]
where the upper bound vanishes to 0 as \(n \to \infty\). To see this, note that the CDF on \(y \in (0,\theta_0]\) is  \(F_Y(y) = y^2/\theta_0\). Then, we have

\vspace{-0.5cm}

\[
P_{\theta_0}(Y_{(1)}> \varepsilon) = (1 - F_Y(\varepsilon))^n = \left(1 - \frac{\varepsilon^2}{\theta_0}\right)^n \le e^{- n \varepsilon^2 / \theta_0}.
\]

\vspace{-0.5cm}
It follows that, for $\theta_0$ fixed, $\min_i y_i^2=O_p{(n^{-1})}$, and similarly, $\min_i (1-y_i)^2=O_p{(n^{-1})}$. 

\begin{figure}[H]
        \centering
    \includegraphics[width=0.95\linewidth]{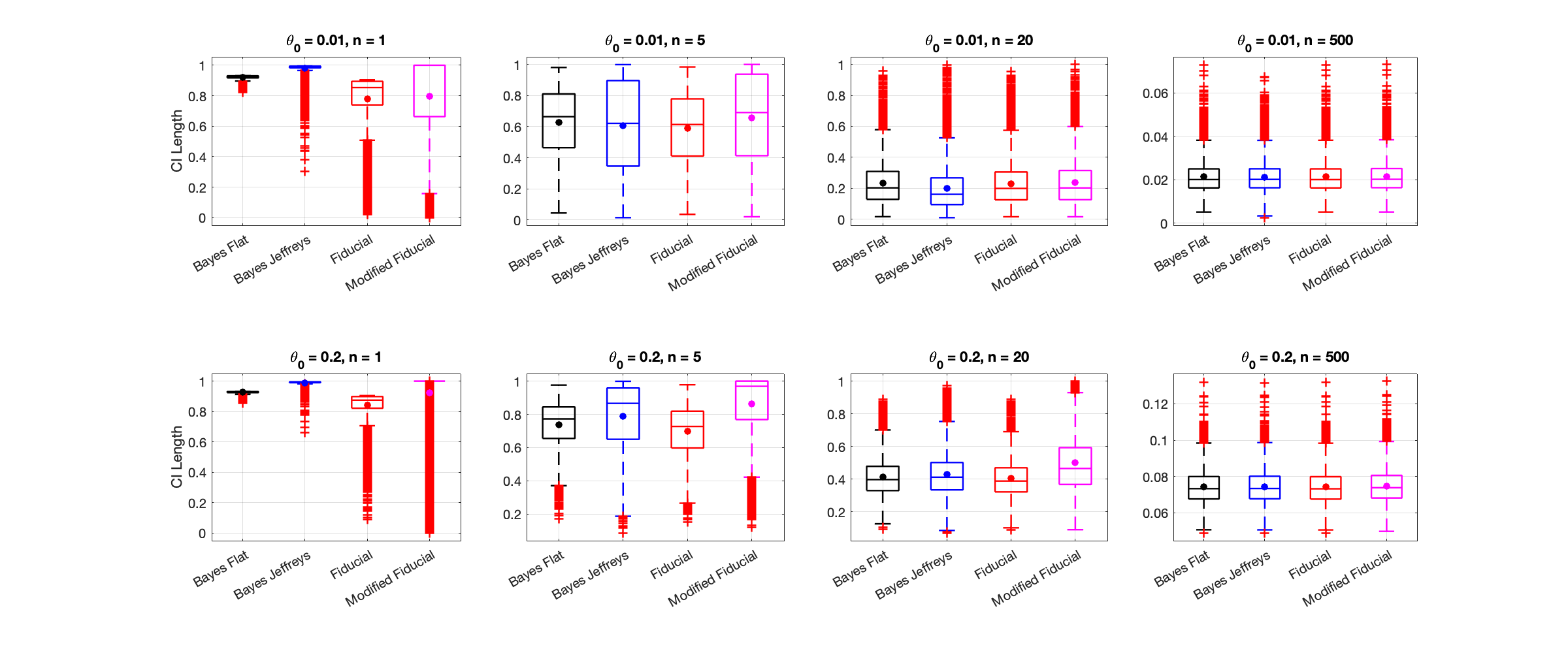}
    \caption{Box-and-whisker plots of two-sided interval lengths for each method at various sample sizes, shown for a true parameter value near the boundary (top row) and in the interior (bottom row).}
    \label{fig:triangle_interval_lengths}
    \end{figure}
\end{example} 

\vspace{-0.5cm}

The blue and magenta lines in Figure~\ref{fig:triangle_coverage} show coverage probabilities corresponding to the Jeffreys Bayesian and modified GFD methods, respectively, which give improved coverage over their alternative counterparts (Bayesian with flat prior and standard GFD) for parameter values near the boundary. Note that the modified GFD achieves exact coverage at $n=1$ (for all parameter values), as expected, while the Jeffreys Bayesian method is conservative there and, in general, more conservative at small sample sizes. Figure~\ref{fig:triangle_interval_lengths} shows box-and-whisker plots of the lengths of the two-sided generalized fiducial and Bayesian posterior credible intervals across sample sizes, for a true parameter value both near the boundary and in the interior. The standard GFD achieves the shortest intervals in general, and has competitive coverage at all sample sizes for interior parameter values, while the modified GFD achieves the best coverage near the boundary and competitive coverage at small sizes in general, but produces slightly wider intervals. The Bayesian approach with Jeffreys prior is most competitive at medium sample sizes for interior parameter values, but is conservative at small sample sizes. As expected from the asymptotic results, all four methods perform similarly for larger sample sizes. Table~1 in the Supplementary Material \citep{supplement} summarizes the empirical coverage and length of two-sided $95\%$ nominal intervals for each method at additional sample sizes and values of the true parameter.

\begin{example}\label{ex:ar2}

As an example of a dependent data setting, we consider a stationary second-order autoregressive [AR(2)] time series model,
\[
Y_t = \phi_1 Y_{t-1} + \phi_2 Y_{t-2} + \sigma \varepsilon_t, 
\quad \varepsilon_t \stackrel{\text{iid}}{\sim} N(0,1),
\]
where the autoregressive coefficients $\boldsymbol{\phi}= (\phi_1,\phi_2)$ satisfy the stationarity conditions:
\[\phi_1 + \phi_2 < 1, \hspace{0.1cm} \phi_2 - \phi_1 < 1, \text{ and } |\phi_2| < 1.\]

The results presented below extend naturally to any AR($p$) model; however, we restrict our attention to the case $p = 2$ for simplicity and to provide a concrete illustration. 

Following the approach in \citet{box1976analysis}, the likelihood function for this model can be expressed as

\vspace{-0.5cm}

\[
p_{\boldsymbol{\phi},\sigma^2}( y_1, \dots, y_T)
=
h(\boldsymbol \phi,\sigma^2)\ \times (\sigma^2)^{-(T-2)/2} \times
\exp\!\left\{
-\frac{1}{2\sigma^2}
\sum_{t=3}^T (y_t - \phi_1 y_{t-1} - \phi_2 y_{t-2})^2
\right\},
\]

\vspace{-0.25cm}

\[
\text{where } \quad
h(\boldsymbol{\phi},\sigma^2)
=
(\sigma^2)^{-1}
\lvert V^{-1} \rvert^{1/2}
\exp\!\left\{
-\frac{1}{2\sigma^2}
(y_1,y_2)' V^{-1} (y_1,y_2)
\right\}
\]
and
\[
V^{-1}
=
\begin{pmatrix}
1 - \phi_2^2 & -\phi_1(1+\phi_2) \\
-\phi_1(1+\phi_2) & 1 - \phi_2^2
\end{pmatrix}.
\]

Let $\boldsymbol{\theta} = (\phi_1,\phi_2,\sigma^2)$, and let $P^{(T)}_{\boldsymbol{\theta}}$ denote the distribution of the full vector $(Y_1,\dots,Y_T)$ with stationary initial distribution for $(Y_1,Y_2)$. According to \citet{kreiss1987adaptive,kreiss1990testing}, the conditional distribution of $(Y_3,\dots,Y_T)$ given $(Y_1,Y_2)$ is locally asymptotically normal with central sequence
\[
\Delta_{T,\boldsymbol{\theta}} =
\frac{1}{\sqrt{T}}
\begin{pmatrix}
\displaystyle \frac{1}{\sigma^2} \sum_{t=3}^T (y_t - \phi_{1} y_{t-1} - \phi_{2} y_{t-2}) 
\begin{pmatrix} y_{t-1} \\ y_{t-2} \end{pmatrix} \\[0.8em]
\displaystyle \frac{1}{2\sigma^4} \sum_{t=3}^T \big[ (y_t - \phi_{1} y_{t-1} - \phi_{2} y_{t-2})^2 - \sigma^2 \big]
\end{pmatrix},
\]
and asymptotic information matrix
\[
I_{\boldsymbol{\theta}} =
\begin{pmatrix}
\frac{1}{\sigma^2}\begin{pmatrix} \gamma_0 & \gamma_1 \\ \gamma_1 & \gamma_0 \end{pmatrix} & 0 \\
0 & \frac{1}{2\sigma^4}
\end{pmatrix},
\]
where $\gamma_k = \mathrm{Cov}(Y_t,Y_{t-k})$ under stationarity. Note that $h(\boldsymbol{\phi},\sigma^2)$ is $O_{P^{(T)}_{\boldsymbol{\theta}}}(1)$ and is therefore negligible in the LAN expansion of the full likelihood.

A DGA for this model can be defined conditionally as: 
\[
G_1(\boldsymbol{\theta}, \boldsymbol \varepsilon) = y_1 = \sqrt{\frac{1-\phi_2}{1+\phi_2}\frac{\sigma^2}{(1-\phi_2)^2-\phi_1^2}}\varepsilon_1 ,
\qquad
G_2(\boldsymbol{\theta}, \boldsymbol \varepsilon) = y_2 = \frac{\phi_1}{1-\phi_2}y_1 + \sqrt{\frac{\sigma^2}{1-\phi_2^2}}\varepsilon_2,
\]
and for \(t \ge 3\),
\[
 G_t(\boldsymbol{\theta}, \boldsymbol \varepsilon) = y_t = \phi_1 y_{t-1} + \phi_2 y_{t-2} + \sigma \varepsilon_t.
\]

Then, to derive the GFD, we invert the DGA to solve for the auxiliary variable, differentiate the data generating function \(G(\boldsymbol{\varepsilon}, \boldsymbol{\theta})\)
with respect to the parameters, and substitute the inversions for \(\varepsilon_t\). Details of that procedure are provided in Section~4 of the Supplementary Material \citep{supplement}. 

Let \(\mathcal{J}(\mathbf y,\boldsymbol{\theta}) 
= 
\nabla_{\boldsymbol{\theta}}G(\mathbf{u},\boldsymbol{\theta})
\bigg|_{\boldsymbol \varepsilon = G^{-1}(\mathbf y,\boldsymbol{\theta})}\) denote the gradient matrix with columns corresponding to
\(\phi_1, \phi_2, \sigma\) and recall that the Jacobian function in the fiducial solution can be expressed as

\vspace{-0.5cm}

\[
J(\mathbf y,\boldsymbol{\theta}) =
\det\!\left(
\frac{1}{T}
\mathcal{J}(\mathbf y,\boldsymbol{\theta})^\top
\mathcal{J}(\mathbf y,\boldsymbol{\theta})
\right)^{1/2}.
\]

We summarize computation of the Jacobian function here, with full details provided in Section~4 of the Supplementary Material \citep{supplement}. A factor of $1/\sigma$ appearing in the $\sigma$-column of \(\mathcal{J}(\mathbf y,\boldsymbol{\theta})\) can be factored out of the determinant; moreover, this column is a linear combination of the first two columns and
\(
\mathbf{y} = (y_1, y_2, \ldots, y_T)^\top,
\)
so may be replaced by
$\mathbf{y}$ in the determinant. Hence, 
\[
\tilde{\mathcal{J}}(\mathbf y,\boldsymbol{\theta}) =
\begin{pmatrix}
c_1 y_1 & d_1 y_1 & y_1 \\[6pt]
c_2 y_1 & d_2 y_1 + \phi_2 y_2 & y_2 \\[6pt]
y_2 & y_1 & y_3 \\[6pt]
y_3 & y_2 & y_4 \\[6pt]
\vdots & \vdots & \vdots \\[6pt]
y_{T-1} & y_{T-2} & y_T
\end{pmatrix},
\]
where 
\[
c_1 = \frac{\phi_1}{(1-\phi_2)^2-\phi_1^2},
\quad
d_1 = \frac{\phi_1^2 + (1-\phi_2)^2 \phi_2}
{(1-\phi_2^2)\big((1-\phi_2)^2-\phi_1^2\big)}, 
\quad
c_2 = \frac{1}{1-\phi_2},
\text{ and }
d_2 = \frac{\phi_1 \phi_2}{(1-\phi_2)^2}.
\]
This leads to, $J(\mathbf y,\boldsymbol{\theta})
= (\sigma \sqrt{T})^{-1}
\det\!\left(\tilde{\mathcal{J}}^\top \tilde{\mathcal{J}}\right)^{1/2}.$

Section~4 of the Supplementary Material \citep{supplement} contains detailed verification of the conditions for Theorem~\ref{thm:BvM} to show that, under the true parameter $\boldsymbol{\theta_0}$, 
\[
\norm{R_{\sqrt{T}(\bar{\boldsymbol \Theta}-\boldsymbol{\theta}_0)\mid \mathbf{Y}} 
-N\big(I_{\boldsymbol{\theta}_0}^{-1}\Delta_{T,\boldsymbol{\theta}_0}, I_{\boldsymbol{\theta}_0}^{-1}\big)}_{TV} 
\;\;\overset{P_{\boldsymbol{\theta}_0}^{(T)}}{\longrightarrow} 0.
\] 
The key reason these conditions hold despite the dependence between observations is that, as a consequence of stationarity, the autocovariances of observations from an AR(2) process decay exponentially fast. This exponential decay is used to establish a law of large numbers for the terms in the Jacobian function and the existence of uniformly exponentially consistent tests against local alternatives. 

\begin{figure}[h]
    \centering
    \includegraphics[width=0.85\linewidth]{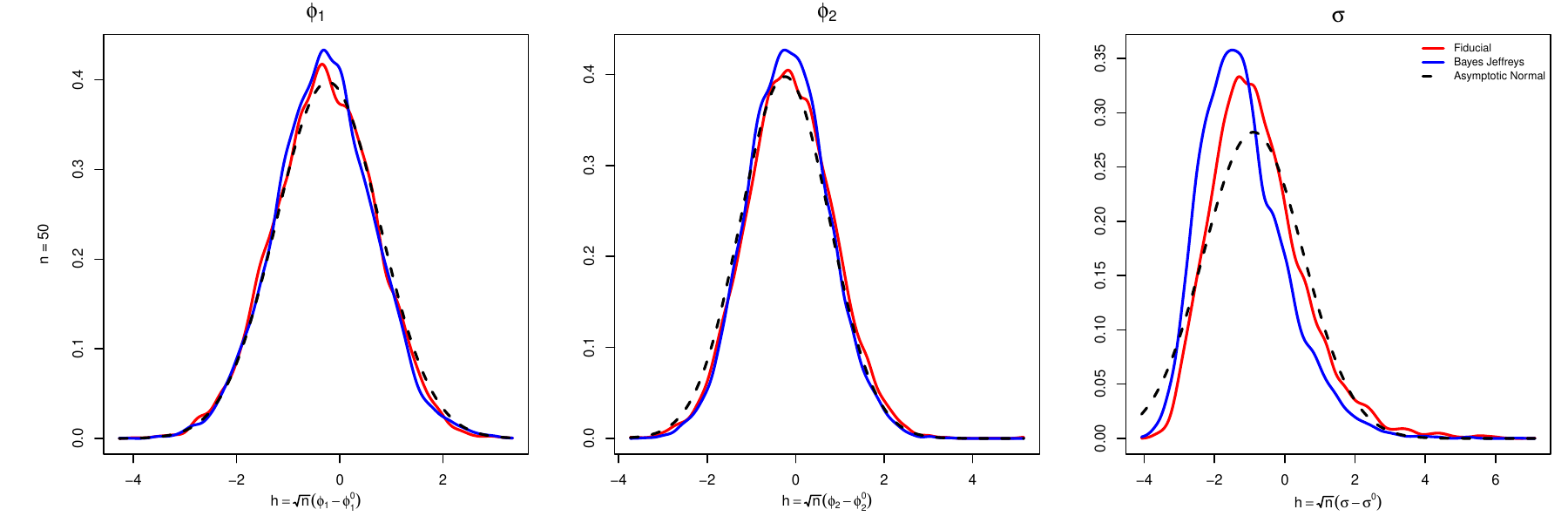}
    \includegraphics[width=0.85\linewidth]{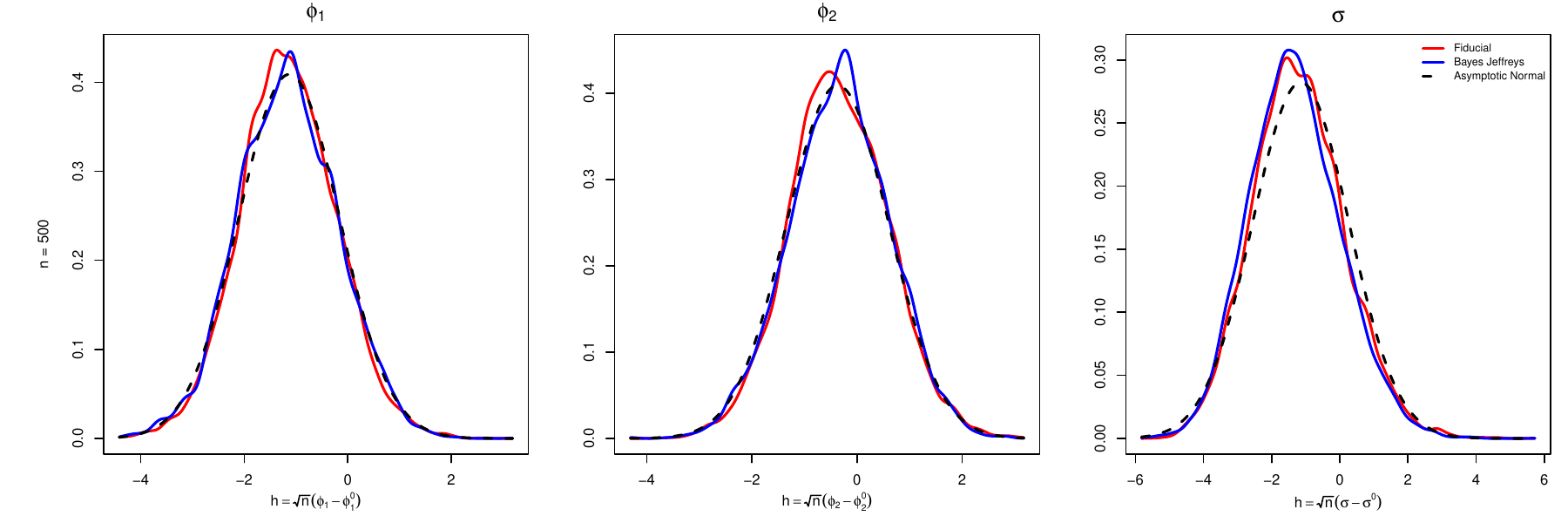}
    \caption{Kernel density estimates for the GFD (red) and Bayesian posterior with Jeffreys prior (blue) for the AR(2) model with true parameters $(\phi_1^0,\phi_2^0,\sigma^0) = (0.5, -0.3, 1)$. Black dashed curves show the corresponding normal approximations $N\big(I_{\boldsymbol{\theta}_0}^{-1}\Delta_{T,\boldsymbol{\theta}_0}, I_{\boldsymbol{\theta}_0}^{-1}\big)$. Top and bottom rows correspond to sample sizes $n=50$ and $n=500$, respectively, with draws rescaled to the local scale; convergence to the normal approximation is evident as $n$ increases from $50$ to $500$.}
    \label{fig:ar2_densities}
\end{figure}

We conducted a simulation study of the AR(2) process, using a Metropolis-within-Gibbs scheme to sample from both the GFD and the Bayesian posterior with Jeffreys prior for parameter estimation. The Fisher information in the Jeffreys prior is computed from the full likelihood including the initial observations; details are provided in Section~4 of the Supplementary Material \citep{supplement}. Figure~\ref{fig:ar2_densities} displays results from the simulation study, where each panel shows kernel density estimates for the GFD (red) and the Bayesian posterior (blue), with the corresponding normal approximations overlaid as black dashed curves. The top row corresponds to a sample size of $n=50$ and the bottom row to $n=500$, with data generated from the true parameter values $\boldsymbol{\theta}_0=(\phi_1^0,\phi_2^0,\sigma^0) = (0.5, -0.3, 1)$. The kernel density estimates are computed after rescaling draws to the local scale, and the normal curve corresponds to $N\big(I_{\boldsymbol{\theta}_0}^{-1}\Delta_{T,\boldsymbol{\theta}_0}, I_{\boldsymbol{\theta}_0}^{-1}\big)$, where $I_{\boldsymbol{\theta}_0}$ is the asymptotic information matrix given earlier. Convergence of the generalized fiducial and Bayesian densities to the normal approximation is observed as the sample size increases from $n=50$ to $n=500$.

Figure~\ref{fig:ar2_phi1_coverage} plots the empirical coverage of nominal $1-\alpha$ two-sided generalized fiducial intervals and Bayesian posterior credible intervals for the AR true parameter $\phi_1^0$ across a range of $\alpha$, based on $M=10,000$ Monte Carlo experiments. The red points and lines correspond to coverage from the GFD, while the blue points and lines correspond to the Bayesian posterior. The black dashed line represents exact coverage. The top panel shows results for a sample of size $n=25$ and the bottom for $n=100$, with data again generated from the true parameters $(\phi_1^0,\phi_2^0,\sigma^0)=(0.5,-0.3,1)$. Figure~\ref{fig:ar2_phi2_coverage} and Figure~\ref{fig:ar2_sigma_coverage} show analogous plots for the parameters $\phi_2^0$  and $\sigma^0$, respectively. For all parameters, Bayesian intervals tend to under-cover at the smaller sample size, but are close to nominal for $n=100$, while intervals from the GFD are closer to nominal in all cases. 

\vspace{-0.8cm}

\begin{figure}[h]
	\centering
	\includegraphics[width=0.8\linewidth]{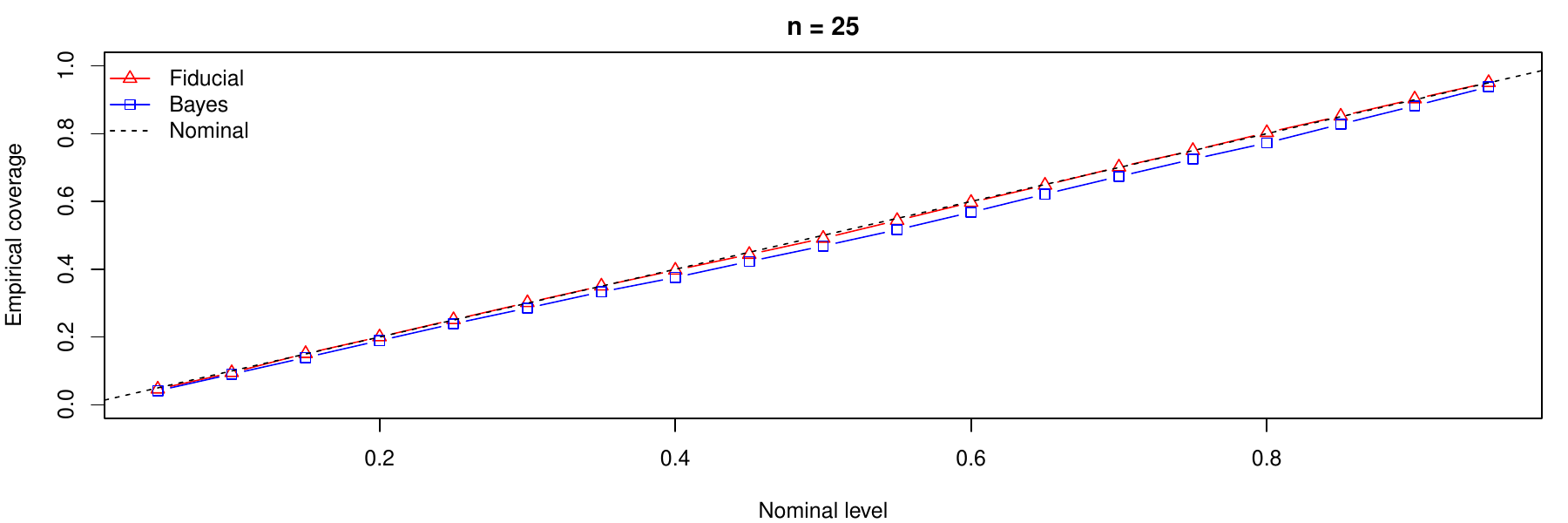}
    \includegraphics[width=0.8\linewidth]{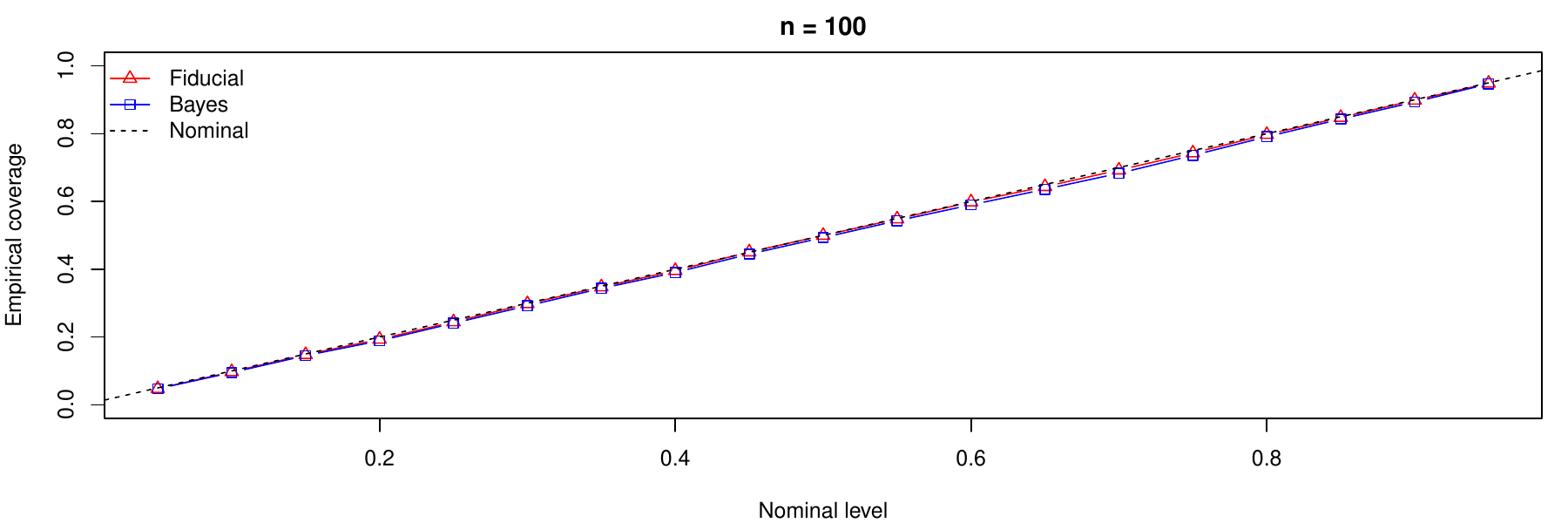}
	\caption{Empirical coverage of two-sided GFD (red) and Bayesian intervals (blue) for the true AR parameter $\phi_1^0$ across nominal levels $1-\alpha$. Exact coverage indicated by a black dashed line. Top and bottom panels correspond to sample sizes $n=25$ and $n=100$, respectively, with data generated from true parameters $(\phi_1^0,\phi_2^0,\sigma^0)=(0.5,-0.3,1)$. Note that Bayesian intervals tend to under-cover compared to the fiducial.}
    \label{fig:ar2_phi1_coverage}
\end{figure}

\vspace{-1cm}

\begin{figure}[h]
	\centering
	\includegraphics[width=0.8\linewidth]{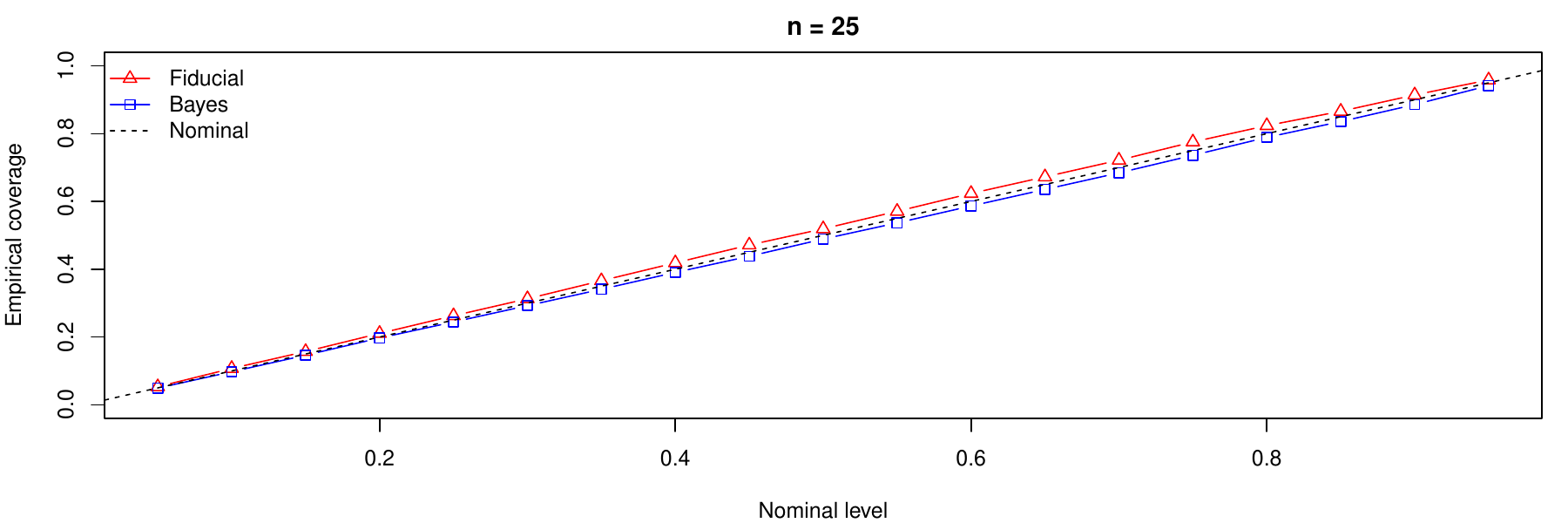}
    \includegraphics[width=0.8\linewidth]{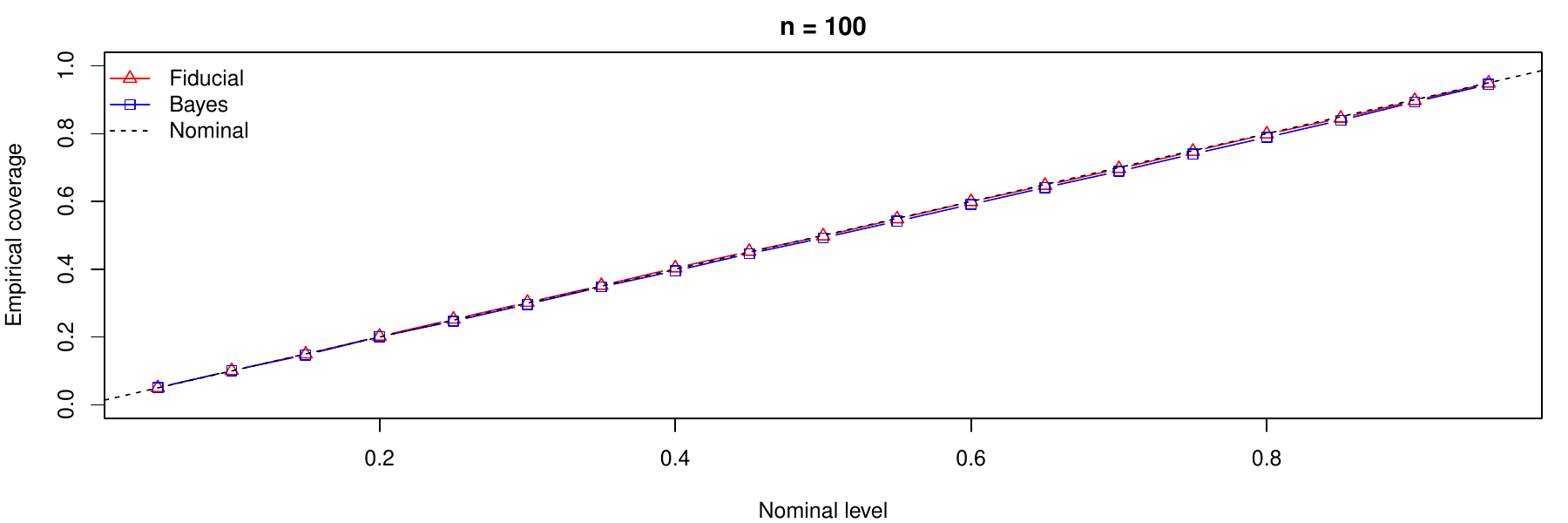}
	\caption{Same as Figure~\ref{fig:ar2_phi1_coverage}, for the AR parameter $\phi_2$.}
    \label{fig:ar2_phi2_coverage}
    \vspace{-5mm}
\end{figure}

\begin{figure}[H]
	\centering
	\includegraphics[width=0.8\linewidth]{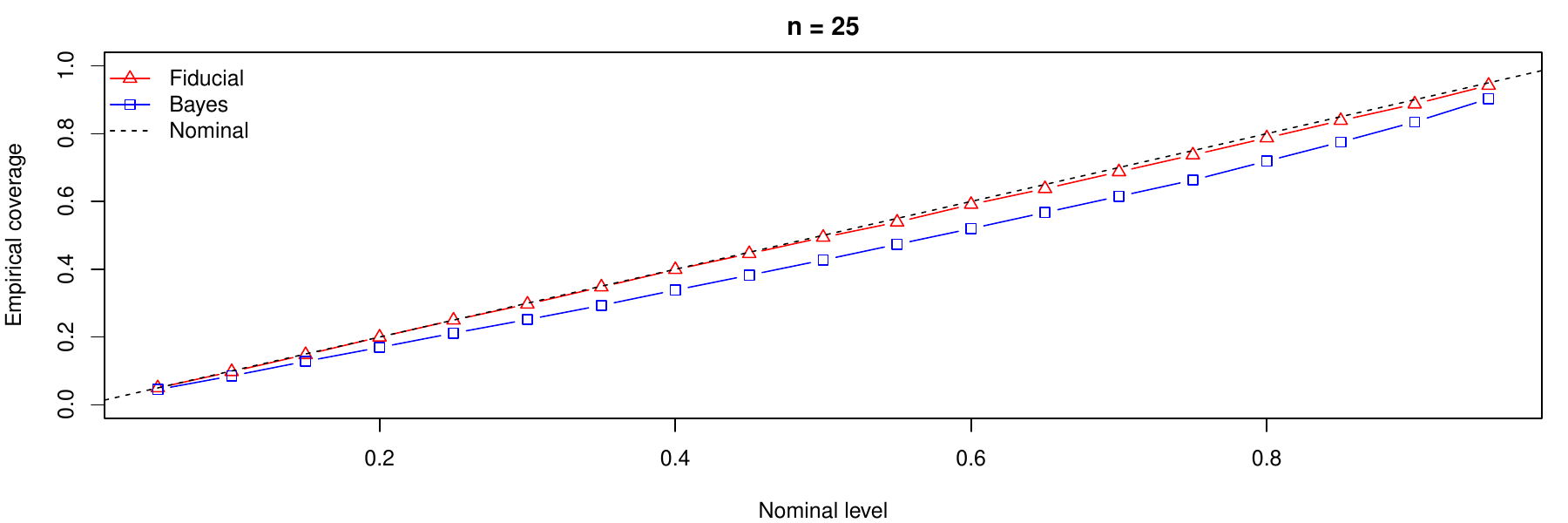}
    \includegraphics[width=0.8\linewidth]{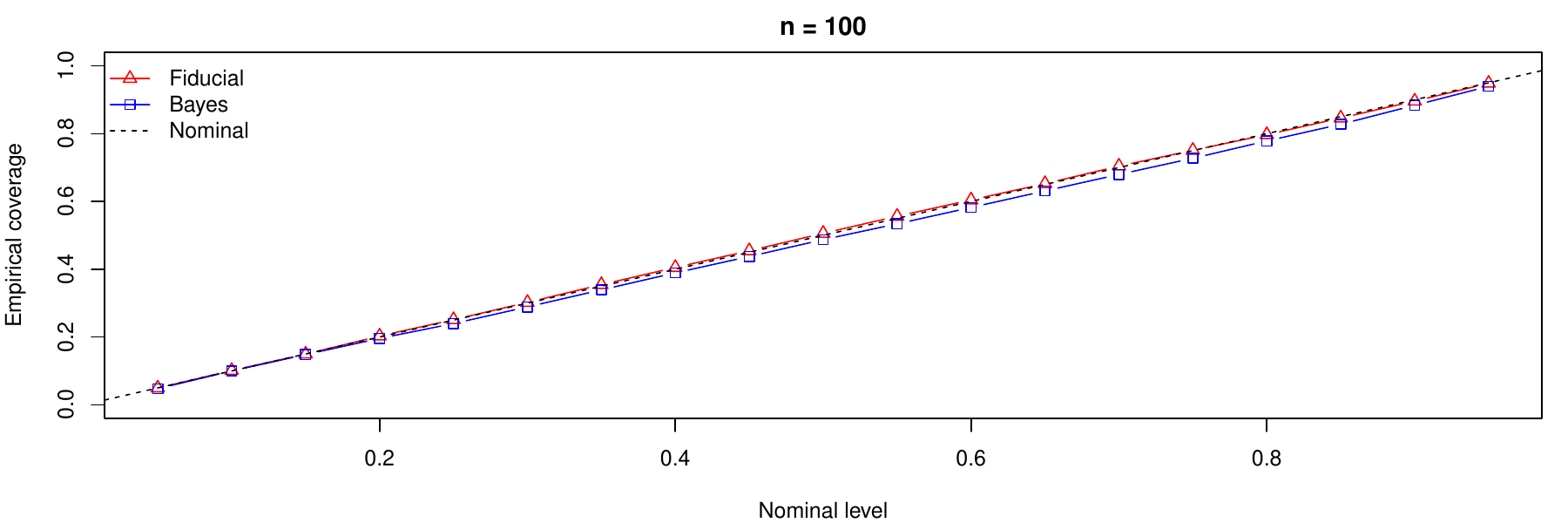}
	\caption{Same as Figure~\ref{fig:ar2_phi1_coverage}, for $\sigma^0$. Under-coverage of Bayesian intervals is most obvious here.}
    \label{fig:ar2_sigma_coverage}
\end{figure}

\end{example}

%
%


\begin{funding}
This research was supported in part by the National Science Foundation under Grant No. DMS-1916115, DMS-2113404, and DMS-2210337.
%
\end{funding}



\bibliographystyle{ba}
\bibliography{references}


\end{document}